\documentclass[12pt,leqno]{amsart}

\usepackage{float}
\usepackage{amssymb,amsmath,amsthm,amscd}
\usepackage{mathrsfs}
\usepackage{enumerate}

\usepackage[all]{xy} 
\CompileMatrices

\numberwithin{equation}{section}

\newcommand{\nc}{\newcommand}
\newcommand{\C}{{\mathbb C}}
\newcommand{\Q}{\mathbb {Q}}

\newcommand{\Z}{{\mathbb Z}}
\newcommand{\B}{{\mathcal{B}}}
\newcommand{\I}{{\mathscr{I}}}

\renewcommand{\L}{{\mathscr{L}}}

\newcommand{\Der}{{\operatorname{Der}}}

\newcommand{\D}{\mathscr{D}}

\newcommand{\Gt}{{\mathfrak{t}}}
\newcommand{\ts}{{\mathfrak{t}}\,}
\newcommand{\gl}{{\mathfrak{gl}}}

\newcommand{\Gr}{{\rm Gr}}

\newcommand{\Pbb}{\mathbb{P}}

\newcommand{\codim}{{\operatorname{codim}}}
\newcommand{\Spec}{{\operatorname{Spec}}}

\newcommand{\seteq}{\mathbin{:=}}

\theoremstyle{plain}
\newtheorem{lemma}{Lemma}[section]
\newtheorem{prop}[lemma]{Proposition}
\newtheorem{theorem}[lemma]{Theorem}
\newcommand{\Prop}{\begin{prop}}
\newcommand{\enprop}{\end{prop}}
\newcommand{\Lemma}{\begin{lemma}}
\newcommand{\enlemma}{\end{lemma}}
\newcommand{\Th}{\begin{theorem}}
\newcommand{\enth}{\end{theorem}}
\newtheorem{corollary}[lemma]{Corollary}
\newcommand{\Cor}{\begin{corollary}}
\newcommand{\encor}{\end{corollary}}
\newtheorem{definition}[lemma]{Definition}
\newtheorem{question}[lemma]{Question}
\newcommand{\Question}{\begin{question}}
\newcommand{\enquestion}{\end{question}}
\newcommand{\Def}{\begin{definition}}
\newcommand{\edf}{\end{definition}}

\theoremstyle{definition}
\newtheorem{remark}[lemma]{Remark}

\nc{\Remark}{\begin{remark}}
\nc{\enremark}{\end{remark}}

\newcommand{\g}{{\mathfrak{g}}}
\newcommand{\Gg}{{\mathfrak{g}}}

\newcommand{\Hom}{\operatorname{Hom}}
\newcommand{\End}{\operatorname{End}}

\newcommand{\isoto}[1][]{\mathop{\xrightarrow[#1]%
{\rule{0pt}{.9ex}%
{\raisebox{-.4ex}[0ex][-.6ex]{$\mspace{3mu}\sim\mspace{3mu}$}}}}}
\renewcommand{\hom}{\operatorname{\it \mathscr{H}\kern-.25em om}}
\newcommand{\tensor}{\mathop\otimes\limits}

\newcommand{\pt}{\operatorname{pt}}

\newcommand{\M}{{\mathscr M}}
\newcommand{\N}{{\mathscr N}}
\newcommand{\eq}{\begin{eqnarray}}
\newcommand{\eneq}{\end{eqnarray}}

\newcommand{\eqn}{\begin{eqnarray*}}
\newcommand{\eneqn}{\end{eqnarray*}}

\newenvironment{tenumerate}{
  \begin{enumerate}
  
  }{\end{enumerate}}

\newcommand{\bnum}{\begin{tenumerate}}
\newcommand{\enum}{\end{tenumerate}}
\newenvironment{anumerate}{
  \begin{enumerate}
  
  }{\end{enumerate}}

\newcommand{\bia}{\begin{anumerate}}
\newcommand{\eia}{\end{anumerate}}
\newcommand{\on}{\operatorname}

\newcommand{\bni}{\begin{tenumerate}}
\newcommand{\eni}{\end{tenumerate}}

\newcommand{\QED}{\end{proof}}
\newcommand{\Proof}{\begin{proof}}
\newcommand{\coh}{{\on{coh}}}

\newcommand{\Supp}{\operatorname{Supp}}
\newcommand{\cl}{\colon}
\newcommand{\To}[1][\phantom{aaaa}]{\xrightarrow{\,#1\,}}
\newcommand{\id}{\on{id}}
\newcommand{\ba}{\begin{array}}
\newcommand{\ea}{\end{array}}

\newcommand{\epi}{\twoheadrightarrow}
\newcommand{\epito}{\twoheadrightarrow}

\newcommand{\set}[2]{\left\{#1 \mathbin; #2 \right\}}

\newcommand{\Mod}{\operatorname{Mod}}
\newcommand{\hs}{\hspace*}

\newcommand{\eqsub}{\begin{subequations}\begin{eqnarray}}
\newcommand{\eneqsub}{\end{eqnarray}\end{subequations}}

\newcommand{\ol}{\overline}

\newcommand{\A}{\mathscr{A}}

\renewcommand{\le}{\leqslant}
\renewcommand{\ge}{\geqslant}

\newcommand{\ext}{{\mathscr{E}xt}}

\newcommand{\SH}{\mathscr{H}}
\newcommand{\SB}{\mathscr{B}}

\nc{\la}{\lambda}
\nc{\lam}{\lambda}
\nc{\U}[1][\g]{U_q(#1)}
\nc{\te}{\tilde{e}}
\nc{\tei}{\tilde{e}_i}
\nc{\tf}{\tilde{f}}
\nc{\tfi}{\tilde{f}_i}
\nc{\tU}{\widetilde U_q(\g)}
\nc{\tE}{\tilde{E}}
\nc{\tF}{\tilde{F}}

\nc{\BZ}{{\mathbb{Z}}}
\nc{\al}{\alpha}
\nc{\qs}{{q}}
\nc{\lan}{\langle}
\nc{\ran}{\rangle}
\nc{\re}{{\mathrm{re}}}
\nc{\wt}{\operatorname{wt}}
\nc{\Uf}[1][\g]{U^-_q(#1)}
\nc{\Ue}{U^+_q(\g)}
\nc{\eps}{\varepsilon}
\nc{\vphi}{\varphi}
\nc{\sphi}{\varphi^*}
\nc{\seps}{\varepsilon^*}

\nc{\nn}{\nonumber}

\nc{\vp}{\varpi}
\nc{\cls}{{\operatorname{cl}}}

\nc{\Wt}{{\operatorname{Wt}}}
\nc{\Us}{U'_q(\g)}
\nc{\La}{\Lambda}
\nc{\ro}{{\rm(}}
\nc{\rf}{{\rm)}}
\nc{\norm}{{\mathrm{norm}}}
\nc{\qbox}{\quad\mbox}
\nc{\braid}{{\mathfrak{B}}}
\nc{\Ad}{\operatorname{Ad}}
\nc{\Aut}{\operatorname{Aut}}
\nc{\dt}[1]{\tilde{\tilde #1}}
\nc{\Sn}{S^{{\mathrm{norm}}}}
\nc{\aff}{{\mathrm{aff}}}
\nc{\rk}{{\mathrm{rk}}}
\nc{\tQ}{\widetilde{Q}}
\nc{\tP}{\widetilde{P}}
\nc{\tW}{\widetilde{\mathscr{W}}}
\nc{\Dyn}{\mathrm{Dyn}}
\nc{\tD}{\widetilde{\Delta}}
\nc{\height}{{\operatorname{ht}}}
\nc{\bl}{\bigl}
\nc{\br}{\bigr}
\nc{\Hecke}{\mathrm{H}}
\nc{\HA}{\Hecke^{\mathrm{A}}}
\nc{\HB}{\Hecke^{\mathrm{B}}}
\nc{\K}{\mathrm{K}}
\newcommand{\scbul}{{\,\raise1pt\hbox{$\scriptscriptstyle\bullet$}\,}}
\nc{\vac}{{\phi}}
\nc{\Bt}{\B_\theta(\g)}
\nc{\be}{\begin{enumerate}}
\nc{\ee}{\end{enumerate}}
\nc{\low}{{\mathrm{low}}}
\nc{\upper}{{\mathrm{up}}}
\nc{\Zodd}{\Z_{\mathrm{odd}}}
\nc{\Ft}[1][n]{\mathbb{P}\mathrm{ol}_{#1}}
\nc{\Ftf}[1][n]{\widetilde{\mathbb{P}\mathrm{ol}}_{#1}}
\nc{\KA}{\on{K}^{\mathrm{A}}}
\nc{\KB}{\on{K}^{\mathrm{B}}}
\nc{\Res}{\on{Res}}
\nc{\Fc}[1][{n,m}]{\mathbf{F}_{#1}}
\nc{\tphi}{\tilde{\varphi}}
\nc{\CO}{\mathscr{O}}
\nc{\disc}{\mathfrak{d}}
\nc{\tr}{\on{tr}}
\nc{\Gb}{\mathfrak{b}}
\nc{\stable}{\mathrm{stable}}
\nc{\X}{\mathfrak{X}}
\nc{\Hilb}{\mathrm{Hilb}}
\nc{\W}{\ensuremath{\mathscr{W}}}
\nc{\Ws}{\ensuremath{\rm W}}
\nc{\opp}{{\on{opp}}}
\newcommand{\prolim}[1][]{\mathop{\varprojlim}\limits_{#1}}
\nc{\corps}{{\mathbf{k}}}
\nc{\h}{\mathrm{\hslash}}
\nc{\fL}[1][{\h}]{\C(\mspace{-1mu}(#1)\mspace{-1mu})}
\nc{\ad}{\mathrm{ad}}
\newcommand{\Endm}{\operatorname{\mathscr{E}\kern-.1pc\mathit{nd}}}
\newcommand{\Endomo}{\operatorname{\mathscr{E}\kern-.1pc\mathit{nd}}}
\nc{\bc}{\bar{\corps}}
\nc{\reg}{{\mathrm{reg}}}
\nc{\ysq}{\mathbf{y}^2}
\nc{\Ch}{\on{Ch}}
\nc{\sketch}{\Proof}
\nc{\Gm}{\mathbb{G}_{\mathrm{m}}}
\nc{\hGm}{\hat{\mathbb{G}}_{\mathrm{m}}}
\nc{\ug}{\widehat{\mathrm{G}}_{\mathrm{m}}}

\nc{\tL}{\widetilde{\mathscr{L}}}
\nc{\Fr}{\mathcal{F}}
\nc{\E}{\mathcal{E}}

\nc{\ord}{\on{ord}}
\nc{\bM}{\overset{\hs{1.5ex}\rule[-.08ex]{1.8ex}{.08ex}}{\M}}
\nc{\romano}{\mathrm{o}}
\nc{\into}{\hookrightarrow}
\nc{\good}{\mathrm{good}}
\nc{\tA}{\widetilde\A}
\nc{\Vz}{{V}\kern-1.1ex\raisebox{1.5ex}[0ex][0ex]{$\cdot$}}
\nc{\bxes}[1]{\raisebox{.9ex}{$\cdot$}{\kern#1}\raisebox{0ex}{$\cdot$}
{\kern#1}\raisebox{-.9ex}{$\cdot$}}
\nc{\ssum}{\mathop{\mbox{\normalsize{${\sum}$}}}\limits}

\usepackage{stmaryrd}
\def\iso{\isoto}
\newcommand{\bbC}{{\mathbb C}}

\newcommand{\Y}{{\mathscr{Y}}}
\def\GL{\operatorname{GL}\nolimits}

\nc{\Fs}{\ensuremath{\rm F}}
\def\can{{\mathrm{can}}}
\nc{\isotf}{\overset{
{\rule{0pt}{.9ex}%
{\raisebox{-.6ex}[0ex][-.7ex]{$\mspace{3mu}\sim\mspace{3mu}$}}}}
{\longleftrightarrow}}
\nc{\tN}{\tilde\N}
\nc{\tens}{\mathop\otimes\limits}

\begin{document}

\title
{Microlocalization of Rational Cherednik algebras
}

\author{Masaki KASHIWARA}
\author{Rapha\"el Rouquier}
\thanks{The first author is partially supported by 
Grant-in-Aid for Scientific Research (B) 18340007,
Japan Society for the Promotion of Science.}

\address{M.K.: Research Institute for Mathematical Sciences,
Kyoto University, Kyoto 606, Japan
}
\address{R.R.:
Mathematical Institute, University of Oxford\\
24-29 St Giles'\\
Oxford, OX1 3LB, UK}

\date{May 8, 2007}

\maketitle


\begin{abstract}
We construct a microlocalization of the rational Cherednik
algebras $H$ of type $S_n$. This is achieved by a quantization of the Hilbert
scheme $\Hilb^n\C^2$ of $n$ points in $\C^2$.
We then prove the equivalence of the category of $H$-modules and the one of
modules over its microlocalization under certain conditions on the parameter.
\end{abstract}

\tableofcontents

\section{Introduction}

\smallskip
Let us recall that 
$\Hilb^n\C^2$, the Hilbert scheme of $n$ points in $\C^2$,
is a symplectic (in particular
crepant) resolution of $\C^{2n}/S_n=S^n\C^2$. On the other hand,
the orbifold $[\C^{2n}/S_n]$ (or the corresponding algebra
$\C[\C^{2n}]\rtimes S_n$) is a non-commutative crepant resolution of
$\C^{2n}/S_n$. There is an equivalence between derived categories
of coherent sheaves on $\Hilb^n\C^2$ and finitely generated modules over
$\C[\C^{2n}]\rtimes S_n$ (McKay's correspondence, cf.\ \cite{Ha}).

The rational Cherednik algebra $H_c$ associated with $S_n$ is a one-parameter
quantization of $\C[\C^{2n}]\rtimes S_n$. We construct a one-parameter
quantization $\tA_c$ of $\CO_{\Hilb^n\C^2}$ and an equivalence of categories
between a certain category of $\tA_c$-modules (good modules with
$F$-action) and the category of finitely generated $H_c$-modules (under
certain conditions on the parameter $c$). Note
that this is an equivalence of abelian categories, while the non-quantized
McKay's correspondence is only an equivalence of derived categories.

The quantization $\tA_c$ is a sheaf over $\Hilb^n\C^2$. Locally on an open
subset isomorphic to $T^*U$, it is isomorphic to the sheaf of 
micro-differential operators $\W$ with a homogenizing parameter
$\hbar$.

\smallskip
Note that our construction is an analog of the Beilinson-Bernstein
localization Theorem for universal enveloping algebras upon flag varieties: 

\medskip
\hs{10ex}\fbox{$$\parbox{50ex}{$\begin{array}{l|l}
\textrm{nilpotent cone } \N & \C^{2n}/S_n \\[1ex]
\textrm{enveloping algebra quotients } U_\lambda(\Gg) & H_c \\[1.5ex]
T^*(G/B) & \Hilb^n\C^2 \\[1ex]
\D_{G/B.\,\lambda} & \tA_c
\end{array}$}$$}

\medskip

$$\xymatrix{
{\D}_{G/B.\,\lambda}\ar@{.>}[d]_{\textrm{quantization}} 
\ar@{<.>}[rr]^-{\sim} &&
U_\lambda(\Gg) \ar@{.>}[d]^{\textrm{quantization}} &&&
{\tA}_c \ar@{.>}[d]_{\textrm{quantization}} \ar@{<.>}[rr]^-{\sim} &&
H_c \ar@{.>}[d]^{\textrm{quantization}} \\
T^*(G/B) \ar[rr]_{\textrm{resolution}} && \N &&&
{\Hilb^n}\C^2 \ar[rr]_{\textrm{resolution}} && \C^{2n}/S_n
}$$

\smallskip
Let us mention that our constructions give rise to the spherical
subalgebra $eH_ce$ of $H_c$ and under certain assumptions on $c$ the
two algebras are Morita equivalent. It would be interesting to
quantize directly the Procesi bundle to obtain $H_c$.

\bigskip
Let us now describe some earlier results related to our work.
An important achievement of Etingof and Ginzburg \cite{EtGi}
and of Gan and Ginzburg \cite{GG} is a construction
of a deformation of the Harish-Chandra morphism for $\GL_n(\C)$,
providing a construction of the spherical subalgebra $eH_ce$ of $H_c$
as a quantum Hamiltonian reduction. This provides
a quantization of the Calogero-Moser space, which is itself obtained by
classical Hamiltonian reduction (Kazhdan, Kostant and Sternberg \cite{KKS}).

Gordon and Stafford \cite{GS,GS2} construct a one-parameter family
of graded ($\BZ$)-algebras $\B_c$ that quantize (a graded ($\Z$)-algebra 
Morita-equivalent to) the homogeneous coordinate
ring of $\Hilb^n\C^2$.

In positive
characteristic, Bezrukavnikov, Finkelberg and Ginzburg \cite{BFG}
construct a sheaf of
Azumaya algebras on the Hilbert scheme whose algebra of global
sections is isomorphic to $H_c$ and obtain an equivalence
of derived categories between modules over that Azumaya algebra and
representations of $H_c$.

\bigskip
Let us explain the type of sheaf of algebras used to quantize
$\Hilb^n\C^2$.
On a complex
contact manifold, Kashiwara \cite{Kc} constructed the stack $\E$ of
microdifferential operators.
Locally, a model for a contact
manifold is the projectivized cotangent bundle
$P^*X$ and the stack $\E$ comes from the sheaf $\E_X$ of
microdifferential operators of Sato, Kawai and Kashiwara.

On a symplectic variety, Kontsevich \cite{Ko} and
Polesello-Schapira \cite{PS} defined a stack $\W$ of microdifferential
operators
with a homogenizing parameter $\hbar$ (making all objects modules over
$\fL$).
Locally, a model is $T^*X$ and $\W$ comes from microdifferential
operators on $P^*(X\times\C)$ which do not depend on the extra variable.

For applications to representation theory, these constructions are
unsatisfactory:
\begin{itemize}
\item the first construction ``forgets about the zero-section''
\item the second construction gives ``too large'' objects (defined
over $\fL$ instead of $\C$).
\end{itemize}

\smallskip
To overcome these difficulties,
we consider here symplectic manifolds $X$ with a $\C^\times$-action that 
stabilizes $\C\omega_X$ with a positive weight.
We consider the case where the
stack $\W$ comes from a sheaf of algebras together with a compatible
action of $\C^\times$ and study the corresponding structure, a
``W-algebra with F-action''.
The category of its modules is defined over $\C$,
as the F-action induces a $\C^\times$-
action on $\fL$ whose invariant field is $\C$.

\bigskip
Let us now describe the structure of the paper.

\smallskip
In the first part of this paper (\S\,\ref{sectionW}),
we study
a general setting for the quantization of symplectic manifolds
$X$ with a $\C^\times$-action that 
stabilizes $\C\omega_X$ with a positive weight.
We first review the theory of W-algebras on symplectic manifolds
(\S\,\ref{secreviewW}). In \S\,\ref{secFactions}, we introduce the notion of
``W-algebra with F-action''. An important point of this construction is that
the category of $\W$-modules with F-action 
on a cotangent bundle (for the canonical structure) is equivalent to
the category of modules over the sheaf 
$\D$ of differential operators.
We adapt in \S\,\ref{secequivariance} the study of equivariance and
its twisted version
for the action of a complex Lie group and we explain how to construct
W-algebras with F-action by symplectic reduction in \S\,\ref{secreduction}.
Finally, in \S\,\ref{secaffinity} we provide sufficient conditions to
ensure $\W$-affinity (a counterpart of Beilinson-Bernstein's result for
$\D$-modules).

\smallskip
Section \ref{secH} is devoted to the construction of
$\D$-modules with an action of the rational Cherednik algebra $H_c$
of type $A_{n-1}$
or of its spherical subalgebra $eH_ce$.
This is related to the constructions of
\cite{BFG,GG}. Let $V=\C^n$ and $\g=\gl_n(\C)$.
We construct (\S\,\ref{secconstruction}) a quasi-coherent
$\D_{\g\times V}$-module $\M_c$ together with
an action of $H_c$,
building on the explicit description of the $\D$-module
arising in Springer's correspondence given in \cite{HK}. We construct
a coherent
$\D_{\g\times V}$-submodule $\L_c$ of $\M_c$ that is stable under
the action of the spherical subalgebra of $H_c$ and we construct a
shift operator (\S\,\ref{secspherical}). This is achieved by reduction to
rank $2$.

\smallskip
In \S\,\ref{secHilb} we construct a W-algebra with F-action on $\Hilb^n\C^2$ by
symplectic reduction from the previous constructions.
After recalling some properties of $\Hilb^n\C^2$ in \S\,\ref{secgeo},
we construct in \S\,\ref{subsec:Ana} a W-algebra $\tA_c$
on $\Hilb^n\C^2$ by symplectic
reduction of $\L_c$ for the action of $\GL_n(\C)$.
In \S\,\ref{secaff}, we present our main results: $\tA_c$-affinity of
$\Hilb^n\C^2$, an isomorphism between global sections of
$\tA_c$ and the spherical algebra and an equivalence
between the category of
good $\tA_c$-modules with $F$-action and the one of
finitely generated modules over the spherical algebra.
We also describe similar results for $H_c$.
So, we have obtained a microlocalization of the rational Cherednik
algebras: we have constructed a W-algebra with F-action
over the Hilbert scheme
whose algebra of global sections is isomorphic to $H_c$ and whose
modules are equivalent to representations of $H_c$.
Those results are obtained under certain assumptions on $c$.
We explain in \S\,\ref{sec:quot}
how to view sections of our W-algebras over open subsets of the
Hilbert schemes as appropriate fractions in the Cherednik algebra.
Finally, we describe explicitly the constructions for $n=2$ in 
\S\,\ref{sec:rank2}.

\smallskip
We thank Pierre Schapira for some useful discussions.
The first author thanks Shigeru Mukai for his comments on Hilbert schemes.
The second author thanks Research Institute for Mathematical Sciences,
Kyoto University, for its hospitality.

\section{F-actions on W-algebras}
\label{sectionW}

\subsection{Notations}
By a manifold $M$, we mean a complex manifold, equipped with the classical
topology and $\CO_M$ is the sheaf of holomorphic functions.
We denote by $\D_M$
the sheaf of differential operators with holomorphic coefficients and by
$\E_M$ the sheaf of formal micro-differential operators on 
the cotangent bundle $T^*M$.

We denote by $\Gm$ the multiplicative group $\C^\times$.

Given a ring $A$, we denote by $\Mod_{\coh}(A)$ the category of coherent
left $A$-modules.

\subsection{W-algebras}
\label{secreviewW}
We shall review some results on W-algebras.
We refer the reader to \cite{PS} (where
the convergent version is studied, while we use the simpler
formal version).

\subsubsection{}
Let $\corps=\fL$ 
be the field of formal Laurent series in an indeterminate
$\h$ and let $\corps(0)=\C[[\h]]$.
Given $m\in\Z$, we define $\W_{T^*\C^n}(m)$ as the sheaf of
formal series
$\sum_{k\ge -m}\h^{k}a_k$ ($a_k\in \CO_{T^*\C^n}$) 
on the cotangent bundle $T^*\C^n$ of $\C^n$ and we
set $\W_{T^*\C^n}=\cup_{m}\W_{T^*\C^n}(m)$.
Then, $\W_{T^*\C^n}$ has a structure of $\corps$-algebra given by
$$a\circ b=
\sum_{\alpha\in\Z_{\ge0}^n}\h^{\vert\alpha\vert}\dfrac{1}{\alpha!}
\partial_\xi^\alpha a\cdot\partial_x^\alpha b.$$
We have a ring homomorphism
$\D_{\C^n}(\C^n)\to\W_{T^*\C^n}(T^*\C^n)$
given by
$x_i\mapsto x_i$, $\dfrac{\partial}{\partial x_i}\mapsto \h^{-1}\xi_i$.

\subsubsection{}
\label{sectionpropertyWalgebras}
Let $X$ be a complex symplectic manifold with symplectic form $\omega_X$.
We denote by $X^\opp$ the symplectic manifold
$X$ with symplectic form $-\omega_X$.

\smallskip
A {\em \Ws-algebra} is a $\corps$-algebra $\W$
on $X$ such that
for any point $x\in X$, there are
an open neighbourhood $U$ of $x$,
a symplectic map $f\cl U\to T^*\C^n$
and a $\corps$-algebra isomorphism $g\cl\W\vert_U\iso f^{-1}\W_{T^*\C^n}$.

\medskip
A W-algebra $\W$ satisfies the following properties.
\be[(i)]
\item 
The algebra $\W$ is a coherent and noetherian algebra.
\item
$\W$ contains a canonical subalgebra $\W(0)$
which is locally isomorphic to $\W_{T^*\C^n}(0)$ (via the maps $g$).
We set $\W(m)=\h^{-m}\W(0)$.
\item
We have a canonical $\C$-algebra isomorphism
$\W(0)/\W(-1)\isoto \CO_X$ (coming from the canonical isomorphism via the maps
$g$).
The corresponding morphism
$\sigma_m\cl \W(m)\to \h^{-m}\CO_X$ is called the {\em symbol map}.
\item We have
$$\sigma_{0}(\h^{-1}[a,b])=\{\sigma_0(a),\sigma_0(b)\}$$
for any $a$, $b\in\W(0)$. Here $\{\scbul,\scbul\}$
is the Poisson bracket.
\item
The canonical map 
$\W(0)\to\prolim[{m\to\infty}]\W(0)/\W(-m)$
is an isomorphism.
\item
A section $a$ of $\W(0)$ is invertible in $\W(0)$ if and only if
$\sigma_0(a)$ is invertible in $\CO_X$.
\item
Given $\phi$ a $\corps$-algebra automorphism of $\W$,
we can find locally an invertible section $a$ of $\W(0)$
such that $\phi=\Ad(a)$.
Moreover $a$ is unique up to a scalar multiple.
In other words, we have canonical isomorphisms
$$\xymatrix{
\W(0)^\times/\corps(0)^\times\ar[r]^\sim_{\Ad}\ar[d]_\sim
&\Aut(\W(0))\ar[d]^\sim\\
\W^\times/\corps^\times\ar[r]^\sim_{\Ad}&\Aut(\W).
}$$
\item
Let $v$ be a $\corps$-linear filtration-preserving
derivation of $\W$.
Then there exists locally a  section $a$ of $\W(1)$
such that $v=\ad(a)$.
Moreover $a$ is unique up to a scalar.
In other words, we have an isomorphism
$$
\W(1)/\h^{-1}\corps(0)\isoto[\ad]\Der_{\mathrm{filtered}}(\W).
$$
\item
If $\W$ is a W-algebra, then its opposite ring
$\W^\opp$ is a W-algebra on $X^\opp$.
\ee
\smallskip
Conjecturally, (iii), (iv) and (v) characterize $\W(0)$.

Note that two W-algebras on $X$ are locally isomorphic.
\subsubsection{}
Assume there exist $a_i$, $b_i\in\W(0)$ ($i=1,\ldots,n$)
such that
$[a_i,a_j]=[b_i,b_j]=0$
and $[b_i,a_j]=\h\delta_{ij}$. 
They induce a symplectic map
$$f=(\sigma_0(a_1),\ldots,\sigma_0(a_n) ; \sigma_0(b_1),\ldots,\sigma_0(b_n)):
 X\to T^*\C^n.$$
Then, there exists a unique isomorphism
$$\W\iso f^{-1}\W_{T^*\C^n},\ 
a_i\mapsto x_i,\ b_i\mapsto \xi_i.$$
We call $(a_1,\ldots,a_n;b_1,\ldots,b_n)$
{\em quantized symplectic coordinates} of $\W$.

\medskip
Let $M$ be a complex manifold $M$ and $\pi_M\cl T^*M\to M$ the projection.
We can associate canonically a W-algebra
$\W_{T^*M}$ with a morphism
$\pi_M^{-1}\D_M\to\W_{T^*M}$ such that
$$\xymatrix{
{\pi_M^{-1}F_m(\D_M)}\ar[r]\ar[d]_{\sigma_m}&{\W_{T^*M}(m)}\ar[d]^{\sigma_m}\\
{\CO_{T^*M}}\ar[r]_{\h^{-m}}&{\h^{-m}\CO_{T^*M}}
}$$
commutes.
Here, $F(\D_M)$ is the order filtration of $\D_M$.
Note that $\pi_M^{-1}\D_M\to\W_{T^*M}$ decomposes into
$\pi_M^{-1}\D_M\to\E_M\to\W_{T^*M}$.
The ring $\W_{T^*M}$ is flat over $\pi_M^{-1}\D_M$
and faithfully flat over $\E_M$.
In particular, for a coherent $\D_M$-module $\M$,
the characteristic variety $\Ch(\M)$ coincides with
$\Supp(\W_{T^*M}\otimes_{\pi_M^{-1}\D_M}\pi_M^{-1}\M)$.

\medskip
Let $X$ and $Y$ be two symplectic manifolds. The product
$X\times Y$ is also a symplectic manifold.
For a W-algebra $\W_X$ on $X$ and a W-algebra $\W_Y$ on $Y$,
there is a W-algebra $\W_X\boxtimes\W_Y$
on $X\times Y$. Letting $p_1\cl X\times Y\to X$
and $p_2\cl X\times Y\to Y$ be the projections,
$\W_X\boxtimes\W_Y$ contains $p_1^{-1}\W_X\otimes_\corps p_2^{-1}\W_Y$
as a $\corps$-subalgebra,
and is faithfully flat over it.

\smallskip
For a $\W$-module $\M$,
a $\W(0)$-{\em lattice} is a coherent $\W(0)$-submodule $\N$ of $\M$
such that the canonical map
$\W\otimes_{\W(0)}\N\to\M$ is an isomorphism.

We say that a $\W$-module $\M$ is {\em good}
if for any relatively compact open subset $U$ of $X$, there exists a
coherent $\W(0)\vert_U$-lattice of $\M\vert_U$.
The full subcategory of good $\W$-modules is an abelian subcategory of
the category of $\W$-modules.

\medskip
The following fact will be used in this paper
(see \cite[Theorem 1.2.2]{KK}, where the convergent version is proved).
\Lemma\label{lem:vanloc}
Let $r$ be an integer and let $\M$ be a coherent $\W$-module such that
$\ext_{\W}^j(\M,\W)=0$ for any $j>r$.
Then $\SH_S^j(\M)=0$ for any closed analytic subset $S$
and any $j<\codim S-r$.
\enlemma

Let $\bc\seteq\cup_{n>0}\fL[\h^{1/n}]$
be an algebraic closure of $\corps$.
We will sometimes need to replace $\W$
with $\corps'\otimes_\corps\W$
for some field $\corps'$ with $\corps\subset\corps'\subset \bc$.

\subsection{F-actions}
\label{secFactions}
\subsubsection{}
Let $X$ be a symplectic manifold. Consider an action of $\Gm$ on $X$,
viewed as a manifold: $\C^\times\ni t\mapsto T_t\in\Aut(X)$. We assume
$\Gm$ stabilizes the line $\C \omega_X\subset H^0(X,\Omega^2_X)$
with a positive weight $m$, i.e.,
$T_t^*\omega_X=t^m\omega_X$ for all $t\in\C^\times$.

We denote by $v$ the vector field given by the $\Gm$-action:
$v(a)(x)=\frac{d}{dt}a(T_t(x))|_{t=1}$. The Poisson bracket
$\{\scbul,\scbul\}$ is homogeneous of degree $-m$:
$$T_t^*\{a,b\}=t^{-m}\{T_t^*a,T_t^*b\}
\textrm{ and }
v\{a,b\}=\{v(a),b\}+\{a,v(b)\}-m\{a,b\}  \textrm{ for }a,b\in\CO_X.$$

\medskip
Let $\W$ be a W-algebra.

\begin{definition}
An {\em \Fs-action with exponent $m$}
on $\W$ is 
an action
of $\Gm$ on the $\C$-algebra $\W$,
$\Fr_t\cl T_t^{-1}\W\isoto \W$ for $t\in \C^\times$,
such that $\Fr_t(\h)=t^m\h$ and
$\Fr_t(a)$ depends holomorphically on $t$ for any $a\in\W$.
\end{definition}

Let us fix an F-action with exponent $m$ on $\W$.
The $\Gm$-action induces
an order-preserving derivation
$v_F$ of $\W$ given by
$v_F(a)=\frac{d}{dt}\Fr_t(a)\vert_{t=1}$. It satisfies the following
properties:
\eq
\ba{l}
v_F(\h)=m\h,\\[1ex]
\sigma_0(v_F(a))=v(\sigma_0(a))\quad\text{for $a\in\W(0)$.}
\ea\label{eq:Fv}
\eneq

\begin{remark}
Here, F stands for ``Frobenius''.
Note that $v_F$ determines the F-action on $\W$.
However, for a given $v_F$ satisfying \eqref{eq:Fv},
we cannot always find an F-action on $\W$.
\end{remark}

The action of $\Gm$ on $\W$ extends to an action on
$\W[\h^{1/m}]=\corps(\h^{1/m})\otimes_\corps\W$ given
by $\Fr_t(\h^{1/m})=t\,\h^{1/m}$.

\begin{definition}
A $\W[\h^{1/m}]$-module with an \Fs-action 
\ro or simply a $(\W[\h^{1/m}],\Fr)$-module\rf\ 
is a $\Gm$-equivariant $\W[\h^{1/m}]$-module: we have isomorphisms
$\Fr_t\cl T_t^{-1}\M\isoto \M$ for $t\in\C^\times$ and we assume that
\be[{\rm(a)}]
\item
$\Fr_t(u)$ depends holomorphically on $t$  for any $u\in\M$
\ro i.e., there  exist locally 
finitely many $u_i$ such that $\Fr_t(u)=\sum_ia_i(t)u_i$
where $a_i(t)\in\W[\h^{1/m}]$ depends holomorphically on $t$\rf,
\item
$\Fr_t(au)=\Fr_t(a)\Fr_t(u)$ for $a\in\W[\h^{1/m}]$, $u\in\M$,
\item 
$\Fr_t\circ\Fr_{t'}=\Fr_{tt'}$ for $t,t'\in\C^\times$.
\ee
\end{definition}

We denote by $\Mod_F(\W[\h^{1/m}])$ the category of
$(\W[\h^{1/m}],\Fr)$-modules: morphisms are morphisms of
$\W[\h^{1/m}]$-modules compatible with the $\Gm$-action.
We denote by $\Mod_F^\good(\W[\h^{1/m}])$ its full subcategory
of good $(\W[\h^{1/m}],\Fr)$-modules.
These are $\C$-linear abelian categories. Note that if there is a
relatively compact open subset $U$ of $X$ such that
$\C^\times\cdot U=X$, then a good $(\W[\h^{1/m}],\Fr)$-module admits
a coherent $(\W(0)[\h^{1/m}],\Fr)$-lattice.

\smallskip
Assume $X=\{\pt\}$, so that $\W=\corps$. We have an equivalence
$\Mod_F(\W[\h^{1/m}])\iso\Mod(\C)$, $\M\mapsto \M^{\Gm}$, with 
quasi-inverse $V\mapsto \fL[\hbar^{1/m}]\otimes_\C V$.

\begin{remark}
Kontsevich and Kaledin \cite{Kal} have also studied quantization
for a symplectic variety with a $\Gm$-action that stabilizes $\C\omega_X$
with a positive weight.
\end{remark}

\subsubsection{}
Let $\W$ be a W-algebra with an F-action with exponent $m$.
Let $n$ be a positive integer and consider the restriction of the
F-action via $\Gm\to\Gm,\ t\mapsto t^n$: we have a new action given
by $T'_t=T_{t^n}$ and $\Fr'_t=\Fr_{t^n}$.
This defines an F-action on $\W$ with exponent $mn$.
Then, we have quasi-inverse equivalences of categories
\eqn
\Mod_F(\W[\h^{1/m}])&\isotf& \Mod_F(\W[\h^{1/nm}])\\
\M&\mapsto& \W[\h^{1/nm}] \otimes_{\W[\h^{1/m}]}\M \\
\set{s\in\N}{\text{$\Fr'_\zeta(s)=s$ for any $\zeta\in\C$ with
$\zeta^n=1$}}&\mapsfrom & \N
\eneqn

\begin{remark}
The equivalence above shows the category depends only on the $1$-parameter
subgroup of $\Aut(X,\W)$ given by the $\Gm$-action. 

Let $\hGm=\prolim[n] \Gm$, where the limit is taken over maps
$f_{n,n'}\cl\Gm\to\Gm,\ t\mapsto t^{n/n'}$ for positive integers $n,n'$
with $n'|n$. This is a pro-algebraic group (some sort of universal
covering group of $\Gm$). In terms of functions, we have
$\hGm=\Spec(\bigoplus_{a\in\Q}\bbC t^a)$ with multiplication
coming from the coproduct $t^a\mapsto t^a\otimes t^a$. Instead of
considering $\Gm$-actions as above, we could
consider $\hGm$-actions on $X$ such that
$T_t^*\omega_X=t\omega_X$. Although theoretically more satisfactory, 
this more complicated formulation is not used in the present paper.
\end{remark}

\subsubsection{}
Let us now give two examples.

\smallskip
Let $M$ be a manifold,
$X=T^*M$ and $\W=\W_{T^*M}$.
We consider the canonical $\Gm$-action given by
$T_t(x,\xi)=(x,t\xi)$. There is a unique F-action with exponent $1$ on
$\W$ with $\Fr\vert_{\D_M}=\id$.
Then, for any $\Gm$-invariant open subset $U$ of $X$,
we have an equivalence
$$
\Mod_F^\good(\W\vert_U)\iso
\Mod_\good(\E_M\vert_U),\ \M\mapsto \M^{\Gm}.$$
In particular, we have an equivalence
$$\Mod_F^\good(\W)\iso
\Mod_\good(\D_M).$$

\medskip
Let $X=T^*\C^n$ and $\W=\W_{T^*\C^n}$. Fix $m>1$ and
$l_1,\ldots,l_n\in\{1,\ldots,m-1\}$.
We define a $\Gm$-action by
$T_t((x_i),(\xi_i))=((t^{l_i}x_i), (t^{m-l_i}\xi_i))$.
Then $T_t^*(\omega_X)=t^m\omega_X$.
We define an F-action on $\W$ with exponent $m$
by
$\Fr_t(x_i)=t^{l_i}x_i$,
$\Fr_t(\partial_i)=t^{-l_i}\partial_i$,
and $\Fr_t(\h)=t^m\h$ (note that
the relation $[\partial_i,x_i]=1$ is preserved by $\Fr_t$).
Then,
$$\End_{\Mod_F(\W[\h^{1/m}])}(\W[\h^{1/m}])^\opp
=\C[\h^{-l_i/m}x_i,\h^{l_i/m}\partial_i;i=1,\ldots,n]\subset
\W[\h^{1/m}],$$
which is isomorphic to $\D(\C^n)$.
Moreover, $\Mod_F^{\good}(\W[\h^{1/m}])$
is equivalent to $\Mod_{\coh}(\D(\C^n))$
(see Theorem~\ref{th:Feq} below).

\subsection{Equivariance}
\label{secequivariance}
We shall discuss $G$-equivariance of $\W$
by adapting \cite{K1,K} where the $\D$-module version is studied.

\subsubsection{}
Let $G$ be a complex Lie group
acting on a symplectic manifold $X$.
Given $g\in G$, let $T_g$ be the corresponding 
symplectic automorphism of $X$.
Let $\g$ be the Lie algebra of $G$ and 
assume that a moment map
$\mu_X\cl X\to \g^*$ is given.

A W-algebra with $G$-action is
a W-algebra with an action of $G$:
we have $\corps$-algebra isomorphisms
$\rho_g\cl\W\isoto T_g^{-1}\W$ for $g\in G$
such that for any $a\in\W$, $\rho_g(a)$ depends holomorphically
on $g\in G$.
Moreover we assume that there is a {\em quantized moment map}
$\mu_\W\cl \g\to \W(1)$
such that
\eqn
&&\ba{l}
[\mu_\W(A),a]=\dfrac{d}{dt}\rho_{\exp(tA)}(a)\vert_{t=0},\\[1ex]
\sigma_0(\h\mu_\W(A))=A\circ\mu_X,\\[1ex]
\mu_\W(\Ad(g)A)=\rho_g(\mu_\W(A))
\ea\qquad
\txt{for any $A\in\g$ and $a\in \W$.}
\eneqn
Note that $\mu_\W$ is a Lie algebra homomorphism.

\subsubsection{}
\label{sectwisted}
A quasi-$G$-equivariant $\W$-module is a
$\W$-module $\M$ with an action of $G$:
$$\rho_g\cl \M\isoto T_g^{-1}\M$$
depending holomorphically on $g\in G$
and such that $\rho_g(au)=\rho_g(a)\rho_g(u)$
for $a\in\W$ and $u\in\M$.
Then, we have a Lie algebra homomorphism
$\alpha\cl\g\to \End_\corps(\M)$
given by $\alpha(A)(u)=\frac{d}{dt}\rho_{\exp(tA)}u\vert_{t=0}$
for $A\in\g$ and $u\in\M$. It satisfies
$$\alpha(A)(au)=[\mu_\W(A),a]u+a\cdot\alpha(A)(u).$$
It follows that we have a Lie algebra homomorphism
\eq
\gamma_\M\cl\g\to\End_\W(\M),\ 
A\mapsto\alpha(A)-\mu_\W(A).\label{eq:gamma}
\eneq

The $\W$-module $\W$ is regarded as a 
quasi-$G$-equivariant $\W$-module. We have
$\alpha(A)=\ad(\mu_\W(A))$ and 
$\gamma_\W(A)(a)=-a\mu_\W(A)$ ($a\in \W$, $A\in\g$).
Given a $G$-module $V$ and a quasi-$G$-equivariant $\W$-module
$\M$, the tensor product $\M\otimes V$ has a natural structure of
a quasi-$G$-equivariant $\W$-module.
The corresponding $\gamma$ is given by
$$\gamma_{\M\otimes V}(A)(u\otimes v)=\gamma_{\M}(A)u\otimes v+u\otimes Av
\quad \text{for $u\in \M$, $v\in V$ and $A\in\g$.}$$
Let $\lambda\in(\g^*)^G$. If
$\gamma_\M$ coincides with the composition 
$\g\To[\lambda]\C\xrightarrow{z\mapsto z\cdot\mathrm{Id}_\M} \End_\W(\M)$,
we say that $\M$ is a twisted $G$-equivariant
$\W$-module with twist $\lambda$.
For such a coherent module $\M$, we have $\Supp(\M)\subset\mu_X^{-1}(0)$.

\smallskip
We denote by $\Mod(\W,G)$ the category of quasi-$G$-equivariant $\W$-modules,
and by $\Mod_\la^{G}(\W)$ its full subcategory of
twisted $G$-equivariant $\W$-modules with twist $\lambda$.
We denote by  $\Mod_\la^{G,\,\good}(\W)$
the category of good twisted $G$-equivariant
$\W$-modules with twist $\lambda$.

The embedding $\Mod_\la^G(\W)\to \Mod(\W,G)$  has a left adjoint
\eq
&&\ba{c}
\Phi_\la\cl\Mod(\W,G)\to \Mod_\la^G(\W)\\[2ex]
\Phi_\la(\M)=
\M/\bigl(\ssum_{A\in\g}(\gamma_\M(A)-\lambda(A))\M\bigr).
\ea
\label{eq:Phi}
\eneq

Let $V$ be a one-dimensional $G$-module
and $\chi\in(\g^*)^G$ its infinitesimal character.
Then, we have an equivalence
\eq
&&\Mod_{\la}^G(\W)\iso \Mod_{\la+\chi}^G(\W),
\ \M\mapsto \M\otimes V.
\label{eq:lachi}
\eneq

Let $\W$ be a W-algebra with an F-action with exponent $m$.
A $G$-action on $(\W,\Fr)$ is a $G$-action on
$\W$ such that
$T_t$ and $T_g$ commute,
$\Fr_t$ and $\rho(g)$ commute and
$\mu_\W(A)$ is $\Fr_t$-invariant, for
$t\in\C^\times$, $g\in G$ and $A\in\g$.

We define similarly the notion of twisted $G$-equivariant
$(\W[\h^{1/m}],\Fr)$-modules.
We denote by $\Mod_{F,\,\la}^{G,\,\good}(\W[\h^{1/m}])$ the category
of good twisted $G$-equivariant $(\W[\h^{1/m}],\Fr)$-modules
with twist $\lambda\in (\g^*)^G$.

\subsection{Symplectic reduction}

\label{secreduction}
Let $X$ be a symplectic manifold with
a symplectic action of $G$ and a moment map
$\mu_X\cl X\to\g^*$.
{\em Assume that $G$ acts properly and freely on $X$}
(i.e.,\ the map $G\times X\to X\times X$ defined by $(g,x)\mapsto (gx,x)$
is a closed embedding).
Then, $\mu_X^{-1}(0)$ is an involutive submanifold.
Let $Z=\mu_X^{-1}(0)/G$,
and let $p\cl \mu_X^{-1}(0)\to Z$ be the projection.
Then $Z$ carries a natural symplectic structure
such that $p$ preserves the symplectic form
(i.e., denoting by $\omega_Z$ the symplectic form of $Z$,
we have $p^*\omega_Z=\omega_X\mid_{\mu_X^{-1}(0)}$).
The local form of $X$ is given by the following Lemma
\cite[\S 41]{GuiSte}.

\Lemma
Locally on $Z$, the manifold $X$ is isomorphic to
$T^*G\times Z$.
More precisely, for any point $x\in\mu_X^{-1}(0)$,
there exist a $G$-invariant open neighbourhood $U$ of $x$ in $X$ and
a $G$-equivariant open symplectic embedding
$U\to T^*G\times T^*\C^n$ compatible with the moment maps.
\enlemma

Let $\W$ be a W-algebra on $X$ with a $G$-action.
Let $\la\in (\g^*)^G$.
Set 
$$\L_\la\seteq\Phi_\la(\W)=\W/\ssum_{A\in \g}\W(\mu_\W(A)+\la(A)).$$
Then, $\L_\la$ is a coherent twisted $G$-equivariant $\W$-module
with twist $\la$. 
The support of $\L_\la$ coincides with $\mu_X^{-1}(0)$.
Let $\L_\la(0)$ be the $\W(0)$-lattice
$\W(0)/\ssum_{A\in \g}\W(-1)(\mu_\W(A)+\la(A))$ of $\L_\la$.

Let $\W_Z=\bigl((p_*\Endomo_\W(\L_\la))^G\bigr)^\opp$,
a sheaf of $\corps$-algebras on $Z$.
\Prop \label{prop:symred}
\bnum
\item
$\W_Z$ is a W-algebra on $Z$, and
$\W_Z(0)\simeq \bigl((p_*\Endomo_{\W(0)}(\L_\la(0)))^G\bigr)^\opp$.
\item
We have quasi-inverse equivalences of categories
\eqn
\Mod^\good(\W_Z)&\isotf&\Mod_\la^{G,\,\good}(\W)\\
\N&\mapsto& 
\L_\la\otimes_{p^{-1}\W_Z}p^{-1}\N\\
(p_*\hom_\W(\L_\la,\M))^G & \mapsfrom & \M.
\eneqn
\item Let $V$ be a one-dimensional representation with infinitesimal character $\chi$.
Then $\N_{\la,\chi}(0)\seteq
(p_*\hom_{\W(0)}(\L_\la(0),\L_{\la-\chi}(0)\otimes V))^G$
is a $\W_Z(0)$-lattice of
a coherent $\W_Z$-module $\N_{\la,\chi}\seteq
(p_*\hom_{\W}(\L_\la,\L_{\la-\chi}\otimes V))^G$
and $\N_{\la,\chi}(0)/\h\N_{\la,\chi}(0)$
is isomorphic to 
$\bigl(p_*(\CO_{\mu_X^{-1}(0)}\otimes V)\br)^G$,
the line bundle on $Z$ associated with $V$.
\item
Assume that $\W$ has an F-action with exponent $m$
compatible with the $G$-action.
Then $\W_Z$ has a natural F-action with exponent $m$ and
we have an equivalence of categories:
$$\Mod_F^\good(\W_Z[\h^{1/m}])\simeq\Mod_{F,\,\la}^{G,\good}(\W[\h^{1/m}]).$$
\enum
\enprop
Note that $\hom_\W(\L_\la,\M)\simeq
p^{-1}((p_*\hom_\W(\L_\la,\M))^G)$.
Hence, if $G$ is connected, we have
$p_*\hom_\W(\L_\la,\M)\simeq (p_*\hom_\W(\L_\la,\M))^G$.

\subsection{$\W$-affinity}
\label{secaffinity}
\subsubsection{}
Let $X$ be a symplectic manifold.
Let $S$ be a variety, let
$f\cl X\to S$ be a projective morphism,
and let $L$ be a relatively ample line bundle on $X$.
Let $\W$ be a W-algebra on $X$.
The following theorem is an analogue of
the result of Beilinson-Bernstein \cite{BB} on $\D$-modules on
flag manifolds. We follow the formulation of \cite{K}.
\Th\label{gth:van}
For $n>0$, let $\L_n(0)$ be a locally free $\W(0)$-module of rank $1$
such that $\L_n(0)/\h\L_n(0)=L^{\otimes(-n)}$.
Set $\L_n=\W\otimes_{\W(0)}\L_n(0)$.

Consider the conditions:
\eq
&&\parbox{70ex}{for $n\gg0$, there exists
a vector space $V_n$ and
a split epimorphism
$\L_n\otimes V_n\epi \W$,
i.e., $\W$ is a direct summand of the direct sum of finitely many
copies of $\L_n$;}
\label{cond:van}\\[1ex]
&&\parbox{70ex}{for $n\gg0$, 
there exists
a vector space $V_n$ and
an epimorphism
$\W\otimes V_n\epi \L_n$.}\label{cond:global}
\eneq

\bnum
\item
Assume \eqref{cond:van}.
Then, for every good $\W$-module $\M$, we have
$R^if_*(\M)=0$ for $i\not=0$.
\item
Assume \eqref{cond:global}.
Then, every good $\W$-module
is generated by its global sections
\ro locally on $S$\rf.
\enum
\enth
The proof will be given in the next two subsections.

\medskip
Assume that $\W$ has an F-action with exponent $m$
and that $S$ has a $\Gm$-action such that $f$ is $\Gm$-equivariant.
Assume moreover that there exists
$\romano\in S$ such that every point of $S$
shrinks to $\romano$
(i.e., $\lim\limits_{t\to0}tx=\mathrm{o}$ for any $x\in S$).

Let $\tW=\W[\h^{1/m}]$ and $A=\End_{\Mod_F(\tW)}(\tW)^\opp$.
\Th\label{th:Feq}
Assume Conditions \eqref{cond:van} and \eqref{cond:global} hold.
Then, $A$ is a left noetherian ring and we have quasi-inverse equivalences
of categories between
$\Mod_F^\good(\tW)$ and $\Mod_{\coh}(A)$ 
\eqn
\Mod_F^\good(\tW)& \isotf & \Mod_{\coh}(A) \\
\M&\mapsto& \Hom_{\Mod_F^\good(\tW)}(\tW,\M)\\
\tW\otimes_AM & \mapsfrom & M.
\eneqn
\enth
The proof will be given in \S\,\ref{sec:Feq}.

\subsubsection{Vanishing theorem}
Let $\W$ be a W-algebra on a symplectic manifold $X$.
Let $\M$ be a coherent $\W$-module.
Recall that $\M(0)$ is a {\em $\W(0)$-lattice} of $\M$
if $\M(0)$ is a coherent $\W(0)$-submodule of $\M$
such that $\W\otimes_{\W(0)}\M(0)\isoto \M$.

\smallskip
We start with the following lemma.
\Lemma\label{lem:lim}
For any coherent $\W(0)$-module $\N$,
the canonical map is an isomorphism
\eq
&&\N\isoto\prolim[m]\N/\h^m\N.
\eneq
\enlemma
\sketch
Let us first show that $\N\to \prolim[m]\N/\h^m\N$ is a monomorphism.
For any $x\in X$, we have morphisms of $\W(0)_x$-modules:
$$\N_x\to\bigl(\prolim[m]\N/\h^m\N\bigr)_x
\to \prolim[m]\bigl(\N_x/\h^m\N_x\bigr).$$
Since the composition is injective (Artin-Rees argument, see e.g. \cite{S}),
the map
$\N_x\to\bigl(\prolim[m]\N/\h^m\N\bigr)_x$ is injective.

\smallskip
Let us show now that $\N\to \prolim[m]\N/\h^m\N$ is surjective.
The question being local,
we can take an exact sequence of coherent $\W(0)$-modules
$$0\to\M\to\L\to\N\to0,$$
where $\L$ is a free $\W(0)$-module of finite rank.
For any Stein open subset $U$ and $m>0$,
we have 
$$H^1(U,\M/(\h^m\L\cap\M))=0,$$
and
$$\text{$\Gamma(U,\M/(\h^m\L\cap\M))
\to\Gamma(U,\M/(\h^{m-1}\L\cap\M))$ is surjective.}$$
Indeed, in the exact sequence
\eqn
&&\Gamma(U;\M/(\h^m\L\cap\M))\to\Gamma(U;\M/(\h^{m-1}\L\cap\M))\\
&&\hs{10ex}\to H^1(U;(\h^{m-1}\L\cap\M)/(\h^m\L\cap\M))\\
&&\hs{20ex}\to H^1(U;\M/(\h^m\L\cap\M))\to H^1(U;\M/(\h^{m-1}\L\cap\M)),\eneqn
$H^1(U;(\h^{m-1}\L\cap\M)/(\h^m\L\cap\M))$ vanishes 
because $(\h^{m-1}\L\cap\M)/(\h^m\L\cap\M)$ is a coherent $\CO_X$-module.

It follows that the following sequence is exact
$$0\to\Gamma(U,\M/(\h^m\L\cap\M))\to
\Gamma(U,\L/\h^m\L)\to
\Gamma(U,\N/\h^m\N)\to0.$$
Since $\{\Gamma(U,\M/(\h^m\L\cap\M))\}_m$ satisfies the ML condition,
the bottom row of the following commutative diagram is exact
$$\xymatrix@C=3ex{
&&{\Gamma(U,\L)}\ar[r]\ar[d]^(.4)\sim&
{\Gamma(U,\N)}\ar[d]\\
0\ar[r]&{\raisebox{-3ex}{$\prolim[m]\Gamma(U,\M/(\h^m\L\cap\M))$}}\ar[r]&
{\raisebox{-3ex}{$\prolim[m]\Gamma(U,\L/\h^m\L)$}}\ar[r]&
{\raisebox{-3ex}{$\prolim[m]\Gamma(U,\N/\h^m\N)$}}\ar[r]&0.
}$$
It follows that $\Gamma(U,\N)\to\prolim[m]\Gamma(U,\N/\h^m\N)
\simeq \Gamma(U,\prolim[m]\N/\h^m\N)$
is surjective.
\QED

\Lemma\label{lem:van}
Let $\M$ be a coherent $\W$-module and let $\M(0)$ be
a $\W(0)$-lattice of $\M$.
Set $\M(m)=\h^{-m}\M(0)$ and $\bM=\M(0)/\M(-1)$.
Assume that
$$\text{$H^i(X,\bM)=0$ for $i\not=0$.}$$
Then,
\bnum
\item
the canonical morphism
$$\Gamma(X,\M(0))/\Gamma(X,\M(-m))\To\Gamma(X,\M(0)/\M(-m))$$
is an isomorphism for any $m\ge0$,
\item
$H^i(X,\M(0))=0$ for any $i\not=0$.
\enum
\enlemma
\sketch
Given $m\ge0$, the exact sequence
$$0\to\bM\To[{\;\h^m}]\M(0)/\M(-m-1)\to\M(0)/\M(-m)\to0$$
induces exact sequences
$$\Gamma(X,\M(0)/\M(-m-1))\to \Gamma(X,\M(0)/\M(-m))
\to H^1(X,\bM)$$
and $$H^i(X,\bM)\to H^i(X,\M(0)/\M(-m-1))\to H^i(X,\M(0)/\M(-m)).$$
It follows that $\Gamma(X,\M(0)/\M(-m-1))\to \Gamma(X,\M(0)/\M(-m))$
is surjective for any $m\ge0$
and $H^i(X,\M(0)/\M(-m))=0$ for any $i>0$.
Since $\Gamma(X,\M(0))=\prolim[m]\Gamma(X,\M(0)/\M(-m))$
by Lemma \ref{lem:lim},
we obtain (i).

For $i>0$, we have
$$H^i(X,\M(0))=\prolim[m]H^i(X,\M(0)/\M(-m))=0$$
because $\{H^{i-1}(X,\M(0)/\M(-m))\}_{m}$
satisfies the ML condition.
\QED

\subsubsection{Proof of Theorem \ref{gth:van}}
Let us prove (i).
The question being local on $S$, we may assume that there exists
a $\W(0)$-lattice $\M(0)$ of $\M$.
Set $\bM=\M(0)/\h\M(0)$.
Then, for $m\gg0$, we have $R^if_*(L^{\otimes m}\otimes_{\CO_X}\bM)=0$
for $i\not=0$. It follows that
$$H^i(f^{-1}U,L^{\otimes m}\otimes_{\CO_X}\bM)=0$$
for any $i\not=0$ and any Stein open subset $U$ of $S$. {}From now on,
we assume that $m$ is large enough so that the vanishing above holds.

\smallskip
Let $\A_m=\Endomo_\W(\L_m)^\opp$, a
\Ws-algebra on $X$. We have
$\A_m(0)=\Endomo_{\W(0)}(\L_m(0))^\opp$.
Let $\L_m(0)^*=\hom_{\W(0)}(\L_m(0),\W(0))$, an
$(\A_m(0),\W(0))$-bimodule,
and let $\L_m^*=\hom_{\W}(\L_m,\W)$,
an $(\A_m,\W)$-bimodule.
We have
$$\L_m^*\simeq \A_m\otimes_{\A_m(0)}\L_m(0)^*\simeq
\L_m(0)^*\otimes_{\W(0)}\W.$$
Note that the bimodules $\L_m$ and $\L_m^*$ give inverse
Morita equivalences between $\A_m$ and $\W$.

Let $\M_m(0)=\L_m^*(0)\otimes_{\W(0)}\M(0)$, an
$\A_m(0)$-lattice in the $\A_m$-module
$\M_m=\L_m^*\otimes_{\W}\M$.
We have $\M_m(0)/\h\M_m(0)\simeq L^{\otimes m}\otimes_{\CO_X}\bM$,
hence
$H^i(f^{-1}U,\M_m(0)/\h\M_m(0))=0$ for $i\not=0$.
Lemma \ref{lem:van} (ii)
implies that
$H^i(f^{-1}U,\M_m(0))=0$ for $i\not=0$.
Taking the inductive limit with respect to
Stein open neighbourhoods $U$ of $s\in S$,
we obtain
$H^i(f^{-1}(s),\M_m(0))=0$,
hence
\eq
&&
H^i(f^{-1}(s),\M_m)\simeq\corps\otimes_{\corps(0)}H^i(f^{-1}(s),\M_m(0))=0.
\label{eq:van}
\eneq

By Condition \eqref{cond:van},
$\W$ is a direct summand of a
direct sum of
finitely many copies of the left $\W$-module $\L_m$.
So, $\W$ is a direct summand of a direct sum
of finitely many copies of the right $\W$-module $\L_m^*$ and
$\M$ is a direct summand of a direct sum
of finitely many copies of $\M_m$ (as a sheaf).
Then, \eqref{eq:van} implies that
$H^i(f^{-1}(s),\M)=0$.
This completes the proof of (i).

\medskip
We now prove (ii).
We shall keep the same notations as in the proof of (i).
Since $L$ is relatively ample, given $s\in S$,
there exists a surjective map
$(\CO_X\mid_{f^{-1}(s)})^{\oplus N}\epi (L^{\otimes m}\otimes\bM)
\mid_{f^{-1}(s)}$ for some $N$.
On the other hand, Lemma \ref{lem:van} (i) implies that
$\Gamma(f^{-1}(s),\M_m(0))\to
\Gamma(f^{-1}(s),\M_m(0)/\h\M_m(0))$ is surjective.
Hence we have a morphism 
$\phi_m\cl\A_m(0)^{\oplus N}\vert_{f^{-1}(s)}\to \M_m(0)\vert_{f^{-1}(s)}$
such that the composition 
$\A_m(0)^{\oplus N}\vert_{f^{-1}(s)}\to
(\M_m(0)/\h\M_m(0))\vert_{f^{-1}(s)}$
is an epimorphism.
It follows that $\phi_m$ is an epimorphism.
Thus, there exists an epimorphism
$\A_m^{\oplus N}\vert_{f^{-1}(s)}\epi \M_m\vert_{f^{-1}(s)}$.
By applying the exact functor
$\L_m\otimes_{\A_m}\scbul\cl\Mod(\A_m)\to\Mod(\W)$, we obtain an epimorphism
$\L_m^{\oplus N}\vert_{f^{-1}(s)}\epi \M\vert_{f^{-1}(s)}$.
The assertion follows now from Condition \eqref{cond:global}.

\subsubsection{Proof of Theorem \ref{th:Feq}}
\label{sec:Feq}
By Theorem \ref{gth:van},
$\Mod_F^\good(\tW)\ni\M\mapsto f_*(\M)\in\Mod(f_*(\tW))$ 
is an exact functor.

By the assumption,
$\romano$ has a neighbourhood system
consisting of relatively compact Stein open 
neighbourhoods $U$
such that $U$ is stable by $T_t$ ($0<\vert t\vert\le 1$).
For such an $U$, we have
$S=\bigcup_{t\in\C^*}T_tU$.
For any
$\M\in\Mod_F^\good(\tW)$,
we have
$$\Hom_{\Mod_F(\tW)}(\tW,\M)
=\set{s\in\M(f^{-1}U)}{\text{$s$ is F-invariant}}.$$
Here $s\in\M(f^{-1}U)$ is F-invariant if
$\Fr_t(s)=s$ for any $t\in\C^\times$ with $|t|=1$.

For $s\in \M(f^{-1}U)$,
let 
$$p_n(s)=\dfrac{1}{2\pi\sqrt{-1}}
\int_{\mid t\mid=1}t^{-n}\Fr_t(s)\dfrac{dt}{t}.$$
We have
$s=\sum_np_n(s)$
and $\h^{-n/m}p_n(s)=p_0(\h^{-n/m}s)$ is F-invariant.

\Lemma\label{exact:modF}
$\Hom_{\Mod_F^\good(\tW)}(\tW,\scbul)$
is an exact functor.
\enlemma
\sketch
Let $\vphi\cl\M\to\M'\to0$ be an epimorphism in $\Mod_F^\good(\tW)$
and let $s'\in\M'(f^{-1}U)$ such that $\Fr_t(s')=s'$ for any $t$ with 
$\vert t\vert=1$.
By Theorem \ref{gth:van},
there exists $s\in\M(f^{-1}U)$ such that
$\vphi(s)=s'$.
We have $\vphi(p_0(s))=s'$ and $p_0(s)$ is F-invariant.
\QED

\Lemma\label{lem:Fgen}
 Any $\M\in\Mod_F^\good(\tW)$ is generated 
by F-invariant global sections.
\enlemma
\sketch
By Theorem \ref{gth:van},
$\M$ is generated by global sections $s_i\in\M(f^{-1}U)$.
Then, $\M$ is generated by the $\h^{-n/m}p_n(s_i)$'s.
Indeed, let $\N$ be the submodule of $\M$ generated by
the $p_n(s_i)$'s. This is a coherent submodule of $\M$.
Let $\psi\cl\M\to\M/\N$ be the quotient morphism.
Then $p_n\psi(s_i)=\psi(p_n(s_i))=0$ for any $n$, and hence
$\psi(s_i)=0$. It follows that $\N=\M$.
\QED

We deduce that
$\Hom_{\Mod_F^\good(\tW)}(\tW,\M)$ is an $A$-module of finite
presentation for any $\M\in\Mod_F^\good(\tW)$.

\Lemma $A$ is left noetherian.
\enlemma
\sketch
Let $I$ be a left ideal of $A$.
Let $\I\subset \tW$ be the image
of $\tW\otimes_AI\to\tW$. Note that $\I$ belongs to $\Mod_F^\good(\tW)$.
Since $\tW$ is coherent, there exist finitely many
$a_i\in I$ such that
$\I=\sum \tW a_i$.
We have
$\Hom_{\Mod_F^\good(\tW)}(\tW,\I)=\sum_iAa_i\subset I$
by Lemma~\ref{exact:modF}.
Since we have injective maps
$I\to \Hom_{\Mod_F^\good(\tW)}(\tW,\I)
\into \Hom_{\Mod_F^\good(\tW)}(\tW,\tW)=A$,
we obtain $I=\sum_iAa_i$.
\QED

Since good $(\tW,F)$-modules are generated by F-invariant sections,
$\Hom_{\Mod_F(\tW)}(\tW,\scbul)$ sends
$\Mod_F^\good(\tW)$ to $\Mod_\coh(A)$.

Given $M\in\Mod_\coh(A)$,
the canonical morphism
$$M\to \Hom_{\Mod_F^\good(\tW)}(\tW,\tW\otimes_AM)$$
is an isomorphism because both sides
are right exact functors of $M$ and the morphism is an isomorphism
for $M=A$.

Given $\M\in\Mod_F^\good(\tW)$, the canonical map
$\tW\otimes_A\Hom_{\Mod_F^\good(\tW)}(\tW,\M)\to\M$
is an isomorphism,
because both sides are right exact functors of $\M$ and
$\M$ has a resolution
$\tW^{\oplus m_1}\to\tW^{\oplus m_0}\to\M\to0$
in $\Mod_F^\good(\tW)$ by Lemma \ref{lem:Fgen}.

This completes the proof of Theorem \ref{th:Feq}.

\section{Rational Cherednik algebras and $\D$-modules}
\label{secH}
\subsection{Definitions, notations and recollections}
\subsubsection{}
Let $V=\C^n$, let $G=\GL(V)=\GL_n(\C)$ and let 
$\g=\gl(V)=\gl_n(\C)$.
We denote by $e_{rs}\in\g$ the elementary matrix with $0$ coefficients
everywhere except in row $r$ and column $s$ where the coefficient is $1$.
We denote by $A_{rs}\in\C[\g]$ the corresponding coordinate function.

We denote by $\ts=\C^n$ the Cartan subalgebra of
diagonal matrices of $\g$ and by $W=S_n$ the Weyl group. We denote by
$s_{ij}$ the transposition $(ij)$ for $1\le i\not=j\le n$. We have
$\C[\ts]=\C[x_1,\ldots,x_n]$ and
$\C[\ts^*]=\C[y_1,\ldots,y_n]$.

We put $\disc(x)=\prod_{i<j}(x_i-x_j)\in\C[\ts]$.
We denote by $\g_\reg$ 
the open subset of regular semisimple elements of $\g$ and we put
$\ts_\reg=\ts\cap\g_\reg=\{x\in\ts ; \disc(x)\not=0\}$.

We will identify $\C[\ts]^W$ and $\C[\g]^G$ via the restriction map.

\smallskip
Given $M$ a graded vector space, we denote by $M_k$ its component of degree 
$k$.

\subsubsection{}
\label{subsubnonchar}
Let $X$ be a manifold, $i\cl Y\into X$  a submanifold,
and let $f\cl \M\to\N$ be a morphism of
coherent $\D_X$-modules. Assume
$Y$ is non-characteristic for $\M$ and $\N$
(i.e., for $Z=\Ch(\M)$ or $Z=\Ch(\N)$, we have
$Z\cap T_Y^*X\subset T_X^*X$). 
If $i^*(f)\cl i^*\M\to i^*\N$ is an isomorphism (resp. monomorphism,
epimorphism),
then so is $f$ on a neighbourhood of $Y$ (see e.g.\ \cite[Theorem 4.7]{K2}).

\subsubsection{}\label{sec:delta}
Let $f\in H^0(X;\CO_X)$ be non zero. We denote by
$\delta(f)$ the element $f^{-1}$ of the $\D_X$-module
$\CO_X[f^{-1}]/\CO_X$. So, $\D_X\delta(f)\subset\CO_X[f^{-1}]/\CO_X$.
More generally, let $S$ be a closed subvariety of complete intersection
of codimension $r$ given by
$f_1=\cdots=f_r=0$ for $f_1,\ldots,f_r\in H^0(X;\CO_X)$.
Then $$\SH_S^j(\CO_X)=0 \text{ for } j\not=r \text{ and }
\SH_S^r(\CO_X)\simeq\CO[(f_1\cdots f_r)^{-1}]/\sum_{1\le i\le r}
\CO[(f_1\cdots \hat{f}_i\cdots f_r)^{-1}].$$
We denote the last $\D_X$-module by $\SB_{S\mid X}$.
We denote by $\delta(f_1)\cdots\delta(f_r)$ 
the section $1/(f_1\cdots f_r)$ of $\SB_{S\mid X}$.

\subsection{Construction of some $\D$-modules}
\label{secconstruction}
\subsubsection{}
\label{subsubCherednik}
Given $c\in\C$, we denote by $H_c$ the rational Cherednik algebra
of $(\ts,W)$ with parameter $c$: this is the $\C$-algebra quotient of
$T(\ts^*\oplus\ts)\rtimes W$ by the relations
\eqn
&&[x_i,x_j]=0,\quad [y_i,y_j]=0,\\
&&[y_i,x_j]=cs_{ij}\quad \text{for $i\not=j$,}\\
&&[y_i,x_i]=1-c\sum_{k\not=i}s_{ik}.
\eneqn

\smallskip
We have a vector space decomposition (``PBW-property'')
\cite[Theorem 1.3]{EtGi}
$$H_c=\C[\ts]\otimes\C[\ts^*]\otimes \C[W].$$
There is an injective algebra morphism (given by Dunkl operators)
\cite[Proposition 4.5]{EtGi}
$$\theta_c\cl H_c\hookrightarrow \D(\ts_{\reg})\rtimes W\subset
\End_\C(\C[\ts_\reg])$$
given by the canonical map on $\C[\ts]\rtimes W$ and by
\eq
\theta_c(y_i)=\partial_{x_i}-c\sum_{k\not=i}\dfrac{1}{x_i-x_k}(1-s_{ik}).
\label{eq:Haction}
\eneq
It induces an isomorphism of algebras after localization
$$\C[\ts_\reg]\otimes_{\C[\ts\,]}H_c\isoto \D(\ts_\reg)\rtimes W.$$

We denote by $e\seteq\dfrac{1}{n!}\sum_{w\in W}w\in\C[W]\subset H_c$
and $e_{\det}\seteq\dfrac{1}{n!}\sum_{w\in W}\det(w)w\in\C[W]\subset H_c$
the idempotents corresponding to the trivial representation
and the sign representation of $W$.

\smallskip
We have an injective morphism $\C[\ts]^W\to eH_ce$, $a\mapsto ae$,
and
we identify $\C[\ts]^W$ with its image.
We put $\ysq=\sum_{i=1}^ny_i^2\in H_c$.
Recall that $eH_ce$ is generated by $\C[\ts]^We$ and $\C[\ts^*]^We$ 
(cf.\ e.g. \cite[proof of Proposition 5.4.4]{BFG}).
On
the other hand, we have an isomorphism of $\C[W]$-modules
(cf.\ e.g. \cite[Corollary 4.9]{BEG})
\eq
\bl(\ad(\ysq)\br)^k\cl \C[\ts]_k\iso \C[\ts^*]_k.
\label{eq:ad}
\eneq
It sends $a(x_1,\ldots,x_n)$ to $2^kk!a(y_1,\ldots,y_n)$.
Hence $eH_ce$ is generated by $\C[\ts]^We$ and $\ysq e$.

\smallskip
We denote by $h\mapsto h^*$ the anti-involution of $H_c$  given by
$x_i\mapsto x_i$, $y_i\mapsto -y_i$, $w\mapsto w^{-1}$
($w\in W$).

\subsubsection{}
We will identify $\g$ and $\g^*$ via the $G$-invariant bilinear
symmetric form $\g\times\g\ni (A,A')\mapsto \tr(AA')$.

A pair $(A,z)$ will denote a point of $\g\times V$.
We identify $T^*(\g\times V)$ with $\g\times \g\times V\times V^*$,
and denote accordingly a point in $T^*(\g\times V)$
by $(A,B,z,\zeta)$. Let $\mu\cl T^*(\g\times V)\to\g^*$
be the moment map. It is given by
$\mu(A,B,z,\zeta)=-[A,B]-z\circ\zeta$.

\smallskip
Let us denote by
$$\mu_D\cl\g\to \D_{\g\times V}(\g\times V)$$
the Lie algebra homomorphism associated with the diagonal action of
$G$ on $\g\times V$.
Let us consider the $\D_{\g\times V}$-module
$\L_c=\D_{\g\times V}u_c$ given by the defining equation:
$$(\mu_D(C)+c\tr(C))u_c=0\quad(C\in\g).$$
More formally, we have $\L_c=\D_{\g\times V}/(\D_{\g\times V}
(\mu_D+c\tr)(\g))$ and $u_c$ is the image of $1$ in $\L_c$.

We consider $\L_c$
as a twisted $G$-equivariant $\D_{\g\times V}$-module
with twist $c\tr$,
where $u_c$ is a $G$-invariant section of $\L_c$.
Since any $a\in\C[\g]^G$ commutes with $\mu_D(C)$ ($C\in\g$), the map
$u_c\mapsto au_c$ extends to a $\D_{\g\times V}$-linear endomorphism of
$\L_c$. Hence, $\L_c$ has
a $(\C[\ts]^W\otimes\D_{\g\times V})$-module structure.

The characteristic variety $\Ch(\L_c)$ of $\L_c$ is the
almost commuting variety:
$$\Ch(\L_c)=\mu^{-1}(0)=\set{(A,B,z,\zeta)}{[A,B]+z\circ\zeta=0}.$$
This is a complete intersection in $T^*(\g\times V)$
\cite[Theorem 1.1]{GG}.

\Lemma
\label{vanishcoh}
Let $\g_1$ be the open subset of $\g$ of elements
which have at least $(n-1)$ distinct eigenvalues.
We have
$$\SH^0_{(\g\setminus\g_\reg)\times V}(\L_c)=0 \text{ and }
\SH^1_{(\g\setminus\g_1)\times V}(\L_c)=0.$$
\enlemma

\sketch
Since $\Ch(\L_c)$ is a complete intersection, 
we have (\cite[(2.23)]{K2})
\eq
\text{$\ext^j_{\D_{\g\times V}}(\L_c,\D_{\g\times V})=0$
for $j\not=\codim_{T^*(\g\times V)}\mu^{-1}(0)=n^2$.}
\label{eq:cohLc}
\eneq

Let $\gamma\cl \g\to \Gt/W$ be the canonical map
associating to $A\in\g$ the eigenvalues of $A$.
Let $\tilde\gamma\cl\mu^{-1}(0)\to\Gt/W$ be given by
$(A,B,i,j)\mapsto \gamma(A)$.
Then,
$\tilde\gamma$ is a flat morphism \cite[Corollary 2.7]{GG}.

Let $S$ be a closed subset of $\ts/W$.
Since $\tilde\gamma$ is flat, we have
$$\codim_{T^*(\g\times V)} (\gamma^{-1}(S)\times_{\g}\Ch(\L_c))-
\codim_{T^*(\g\times V)}\Ch(\L_c)=
\codim_{\Gt/W} S.$$

Lemma \ref{lem:vanloc} applied to $\gamma^{-1}(S)\times_{\g}\Ch(\L_c)$ implies
\eqn
\SH^j_{\gamma^{-1}(S)\times V}(\L_c)=0\quad\text{for } j<\codim_{\Gt/W} S
\eneqn
and the lemma follows.
\QED

\subsubsection{}\label{subsec:delta}
Let us recall some constructions and results of \cite{HK}.
Let $\mu_0\cl\g\to \D_{\ts\times\g}(\ts\times\g)$ be the morphism given by
the action of $G$ on $\ts\times\g$: $g\cdot(x,A)=(x,\Ad(g)A)$.
We consider the $\D_{\ts\times\g}$-module
generated by $\delta_0(x,A)$ with the defining equations:
\eqn
&\mu_0(C)\delta_0(x,A)=0\quad\text{for any $C\in\g$,}\\[1ex]
&
\ba{rcl}
(P(A)-P(x))\delta_0(x,A)&=&0\\[1ex]
(P(\partial_A)-P(-\partial_x))\delta_0(x,A)&=&0
\ea\quad\text{for any $P\in \C[\g]^G$.}
\eneqn

Then, $\D_{\ts\times\g}\delta_0(x,A)$ is a simple holonomic
$\D_{\ts\times\g}$-module with support
$\ts\mathop\times_{\ts/W}\g$.
Its characteristic variety
is the set of
$(x,y,A,B)$ such that $[A,B]=0$ and there exists 
$g\in G$ such that 
$\Ad(g)A$ and $\Ad(g)B$ are upper triangular and
$x$ and $y$ are the diagonal components
of $\Ad(g)A$ and $\Ad(g)B$.
Note that $\D_{\ts\times\g}\delta_0(x,A)\subset
\SB_{\ts\times_{\ts/W}\g\vert \ts\times \g}$
by $\delta_0(x,A)\mapsto\prod_{i=1}^n\delta(P_i(x)-P_i(A))$ 
(see \S\,\ref{sec:delta}),
where $P_i\in\C[\g]^G$ ($i=1,\ldots,n$) are the fundamental invariants
given by $\det(1+tA)=\sum_{i=0}^nP_i(A)t^i$.

We will need to consider the $\D_{\ts\times\g\times V}$-module
$\D_{\ts\times\g}\delta_0(x,A)\boxtimes \CO_{V}$,
generated by $\delta(x,A)\seteq\delta_0(x,A)\boxtimes 1$
which satisfies the same equations as $\delta_0(x,A)$
and $\partial_{z_i}\delta(x,A)=0$.
In particular, $\mu_D(C)\delta(x,A)=0$
for any $C\in\g$.

\subsubsection{}
Let us set 
$$q(A,z)=\det(A^{n-1}z,A^{n-2}z,\ldots,Az,z).$$
We have
$q(\Ad(g)A,gz)=\det(g)q(A,z)$ for $g\in G$
and $[\mu_D(C),q(A,z)]=-\tr(C) q(A,z)$ for $C\in\g$.

Consider the $\D_{\ts\times\g\times V}$-module
$\D_{\ts\times\g\times V}q(A,z)^c\delta(x,A)$.
A precise definition is as follows.
Let us consider the left ideal $\mathscr{I}$
of $\D_{\ts\times\g\times V}\otimes\C[s]$ ($s$ being an indeterminate)
consisting of those $P(s)$ such that
$P(m)q(A,z)^m\delta(x,A)=0$ for any $m\in\Z_{\ge0}$.
We now define $\D_{\ts\times\g\times V}q(A,z)^c\delta(x,A)$
as $\bigl(\D_{\ts\times\g\times V}\otimes\C[s]\bigr)
/\bigl(\mathscr{I}+\D_{\ts\times\g\times V}\otimes\C[s](s-c)\bigr)$.
It is a holonomic $\D_{\ts\times\g\times V}$-module.

The element $q(A,z)^c\delta(x,A)$ satisfies
\eqn
&&(\mu_D(C)+c\tr(C))q(A,z)^c\delta(x,A)=0\quad\text{for any $C\in\g$,}\\
&&(P(A)-P(x))q(A,z)^c\delta(x,A)=0\quad\text{for any $P\in\C[\g]^G$.}
\eneqn
We put $v_c=q(A,z)^c\delta(x,A)$.
Let 
$p_0\cl \ts_\reg\times \g\times V \to \g\times V$ be the projection.
Let us consider the $\D_{\g\times V}$-module
$$\M_c=(p_0)_*(\D_{\ts_\reg\times\g\times V}v_c)
=(p_0)_*(\D_{\ts_\reg\times\g\times V}q(A,z)^c\delta(x,A)).$$
By the definition, 
we have an isomorphism $\M_c\isoto j_*j^{-1}\M_c$
where $j\cl\g_\reg\times V\hookrightarrow \g\times V$ is the open embedding.
This is a quasi-coherent $\D_{\g\times V}$-module
whose characteristic variety is contained in 
the almost commuting variety $\mu^{-1}(0)$.

The action of $W$ on $\Gt_\reg$ induces a $W$-action on
$\M_c$. Here, $W$ acts trivially on $v_c$.
Hence, the $\D_{\g\times V}$-module $\M_c$ has a module structure over
$\D(\ts_\reg)\rtimes W$.
Therefore, $H_c$ acts on 
$\M_c$ via the canonical embedding $\theta_c\cl
H_c\into \D(\ts_\reg)\rtimes W$.

\subsection{Spherical constructions and shift}
\label{secspherical}
\subsubsection{}
There is a $\D_{\g\times V}$-linear homomorphism
\eq
&&\iota\cl\L_c\to \M_c,\ u_c\mapsto v_c.
\label{eq:Lcp}
\eneq
We regard $\M_c$
as a twisted $G$-equivariant $\D_{\g\times V}$-module
with twist $c\tr$,
where sections in $\D(\ts_\reg)v_c$ are $G$-invariant.
Then, the morphism above is $G$-equivariant.
Moreover, it is $\C[\ts]^W$-linear.
Hence $\iota$ induces an epimorphism of
$(\D(\ts_\reg)\rtimes W)\otimes\D_{\g\times V}$-modules:
$$\D(\ts_\reg)\otimes_{\C[\ts]^W}\L_c\epito \M_c.$$

\Lemma\label{lem:reg}
The morphism of $\C[W]\otimes\D_{\g\times V}$-modules
$$1\otimes\iota\cl\C[\ts]\otimes_{\C[\ts]^W}\L_c\to\M_c$$ is an isomorphism
on $\g_\reg\times V$.

In particular, the induced morphisms
$\L_c\xrightarrow{u_c\mapsto v_c} e\M_c$
and $\L_c\xrightarrow{u_c\mapsto \disc(x)v_c} e_{\det} \M_c$
are isomorphisms on $\g_\reg\times V$.
\enlemma

\sketch
Let $i\cl\ts_\reg\times V\into \g\times V$ be the embedding. Note that
$i$ is non-characteristic for $\L_c$ and $\M_c$.
Since $G\cdot \ts_\reg=\g_\reg$, it is enough to prove that the canonical
map $\C[\ts_\reg]\otimes_{\C[\ts_\reg]^W}i^*\L_c\to i^*\M_c$ is an isomorphism
(cf.\ \S\,\ref{subsubnonchar}).

We have
$i^*\mu_D(e_{rs})=(A_{rr}-A_{ss})\partial_{A_{rs}}-z_s \partial_{z_r}$.
 It follows that we have an isomorphism
$$\D_{\ts_\reg\times V}/
\bigl(\ssum_{i}\D_{\ts_\reg\times V}(z_i\partial_{z_i}-c)\bigr)\iso i^*\L_c,
\ 1\mapsto i^*u_c.$$

Let $i''\cl\ts_\reg\times \ts_\reg\into \ts\times\g$ be the embedding.
Since the Jacobian
$$\partial(P_1(x),\ldots,P_n(x))/\partial(x_1,\ldots,x_n)$$
is equal to $\disc(x)$ (e.g.\ \cite[Ch.\;V, \S\;5.4, Proposition 5]{Bourbaki}),
we have
an isomorphism

$$i^{\prime\prime*}\D_{\ts\times\g}\delta_0(x,A)
\xrightarrow[\sim]{\delta_0(x,A)\mapsto\sum_w\disc(a)^{-1}\delta(w^{-1}x-a)}
\bigoplus_{w\in W}\D_{\ts_\reg\times\ts_\reg}\delta(w^{-1}x-a)$$
where $\delta(w^{-1}x-a)=\delta(x_{w(1)}-a_1)\cdots\delta(x_{w(n)}-a_n)$.

Let us denote by $i'\cl
\ts_\reg\times \ts_\reg\times V\into \ts_\reg\times\g\times V$ the embedding.
We have an isomorphism
\eq
\label{eq:genM}
i^{\prime *}\D_{\ts_\reg\times\g\times V}v_c
\xrightarrow[\sim]{v_c\mapsto \sum_w v'_w}
\bigoplus_{w\in W}\D_{\ts_\reg\times\ts_\reg\times V}v'_w.
\eneq
where $v'_w=\disc(a)^{c-1}(z_1\cdots z_n)^c\delta(w^{-1}x-a)$ has
the defining equations
\eqn
&&\ba{lll}
&&\bl(\partial_{x_{w(i)}}+\partial_{a_i}
-(c-1)\smash{\sum_{j\not=i}}\dfrac{1}{a_i-a_j}\br)v'_w=0,\\[1ex]
&&(x_{w(i)}-a_i)v'_w=0,\\[1ex]
&&(z_i\partial_{z_i}-c)v'_w=0,
\ea
\qquad\text{for any $i=1,\ldots,n$.}
\eneqn
In particular, we have
\eq
f(x)v'_w=(w^{-1}f)(a)v'_w\quad\text{for any $f\in\C[\ts]$.}
\label{eq:wshift}
\eneq
We obtain finally an isomorphism
$$i^*\M_c
\xrightarrow[\sim]{v_c\mapsto \sum_w v'_w}
\bigoplus_{w\in W}\D_{\ts_\reg\times V}v'_w.$$
This is compatible with the action of $W$, where $w'(v'_w)=v'_{w'w}$. 
Moreover,
each $\D_{\ts_\reg\times V}v'_w$ is isomorphic to $i^*\L_c$
by $v'_w\mapsto u_c$.
Hence we obtain an isomorphism of
$(\D_{\ts_\reg\times V}\otimes\C[W])$-modules
$$i^*\M_c
\isoto \C[W]\otimes i^*\L_c.$$
The composition
$i^*\bl(\C[\ts]\otimes_{\C[\ts]^W}\L_c\br)\to i^*\M_c
\isoto \C[W]\otimes i^*\L_c$
is given by
$a\otimes u_c\mapsto \sum_{w\in W}w\otimes (w^{-1}a)u_c$
in virtue of \eqref{eq:wshift}.
Then the lemma follows from the fact that
$\C[\ts]\otimes_{\C[\ts]^W}\C[\ts_\reg]\to
\C[W]\otimes\C[\ts_\reg]$ 
given by $a\otimes b\mapsto \sum_{w\in W}w\otimes (w^{-1}a)b$
is an isomorphism.
\QED

\Lemma \label{lem:Lc}
The morphism  $\iota\cl\L_c\to\M_c$ is injective and its image is
stable by $eH_ce$. Furthermore, $eH_ce$ acts faithfully on $\L_c$.
\enlemma
\sketch
The injectivity of $\iota$ follows from Lemma \ref{lem:reg},
because $\L_c$ does not have a non-zero submodule
supported in $(\g\setminus\g_\reg)\times V$
by Lemma~\ref{vanishcoh}.

\smallskip
Since
$eH_ce$ is generated by $\C[\ts]^W$ and $\ysq e$ (cf.\ \S\,\ref{subsubCherednik}),
the stability result follows from the following result
(cf.\ \cite[Proposition 5.4.1]{BFG} and \cite[Proposition 6.2]{EtGi}):
\eq
&&\ysq v_c=\Delta_\g v_c.
\label{eq:Deltauc}
\eneq
Here
$\Delta_\g=\sum\limits_{i,j=1,\ldots ,n}\dfrac{\partial^2}{\partial A_{ij}
\partial A_{ji}}$ is the Laplacian on $\g$.

\smallskip
Finally, the faithfulness of the action of $eH_ce$ follows from the
faithfulness of the action of $H_c$ on $H_cv_c\subset\M_c$.
With the notations of the proof of Lemma \ref{lem:reg}, we have an
isomorphism $i^*\M_c\simeq \D_{\ts_\reg\times V}\rtimes W$
compatible with the action of $\D_{\ts_\reg}\rtimes W$,
 and the faithfulness follows
from that of $\theta_c$.
\QED

\Remark
\bnum
\item
In other words, the subalgebra of $\End_{\D_{\g\times V}}(\L_c)$
generated by $\C[\ts]^W$ and by the endomorphism $u_c\mapsto \Delta_\g u_c$
is isomorphic to $eH_ce$.
\item
The action of $eH_ce$ on $\L_c$ can be described as follows.
Let $\kappa_0\cl\C[\ts]^W\isoto\C[\g]^G\into\D(\g)$
and $\kappa_1\cl\C[\ts^*]^W\isoto\C[\g^*]^G\into\D(\g)$ be
the canonical morphisms. 
We have
\eq
&&\ba{l}
(ae) u_c=\kappa_0(a) u_c\quad\text{for $a\in\C[\ts]^W$,}\\[1ex]
(be) u_c=\kappa_1(b^*) u_c\quad\text{for $b\in\C[\ts^*]^W$.}
\ea
\eneq
The first equality is clear.
We have a commutative diagram
\eq
&&
\ba{c}\xymatrix@C=7em{
\C[\ts]_k^W\ar[r]^{\kappa_0}
\ar[d]_{\bl(\ad(\ysq)\br)^k}&\C[\g]^G_k\ar[d]^{\bl(\ad(\Delta_\g)\br)^k}\\
\C[\ts^*]_k^W\ar[r]_{\kappa_1}&\C[\g^*]^G_k.
}\ea
\label{dia:kappa}
\eneq
{}From \eqref{eq:Deltauc} and the first equality, we deduce that
$$\bl(\ad(\Delta_\g)\br)^k(\kappa_0(a))v_c
=(-1)^k \bl(\ad(\ysq)\br)^k(a)v_c$$ for
$a\in \C[\ts]_k^W$. This gives the second equality.
\enum
\enremark

\subsubsection{}
\label{sec:conditions}
The morphism $\iota$ gives rise to an $(H_c\otimes\D_{\g\times V})$-linear
morphism
\eq
&&
H_ce\otimes_{eH_ce}\L_c\to 
\M_c.
\label{eq:Ldelta}
\eneq

Consider the conditions:
\eq
&&H_ceH_c=H_c,
\label{eq:HeH}\\[1ex]
&&
\text{
$eH_ce_{\det}H_ce=eH_ce$ and $e_{\det}H_ceH_ce_{\det}=e_{\det}H_ce_{\det}$.}
\label{eq:HdiscH}
\eneq

\Lemma If \eqref{eq:HeH} is satisfied, then the morphism 
\eqref{eq:Ldelta} is injective.
\enlemma
\sketch
Since  $H_ce$ is a projective $eH_ce$-module,
any coherent submodule of $H_ce\otimes_{eH_ce}\L_c$ 
vanishes as soon as it is zero on
$\g_\reg\times V$ by Lemma~\ref{vanishcoh}.
Hence it is enough to show that the morphism \eqref{eq:Ldelta} is injective on
$\g_\reg\times V$.
Then the result follows from Lemma \ref{lem:reg} and the fact that
the multiplication map gives an isomorphism of right
$(eH_ce\otimes_{\C[\ts]^W}\C[\ts_\reg]^W)$-modules
$$\C[\ts]\otimes_{\C[\ts]^W}{eH_ce}\otimes_{\C[\ts]^W}\C[\ts_\reg]^W\iso
H_ce\otimes_{\C[\ts]^W}\C[\ts_\reg]^W.$$
\QED

\Prop
Condition \eqref{eq:HeH} holds if and only if
$eH_c$ gives a Morita equivalence between
$H_c$ and $eH_ce$. Similarly,
Condition \eqref{eq:HdiscH} holds if and only if
$eH_ce_{\det}$
gives a Morita equivalence between $e_{\det}H_ce_{\det}$
and $eH_ce$.
\enprop

This follows from the following Lemma:
\Lemma\label{lem:Morg}
Let $A$ be a ring, and let $e_1$ and $e_2$ be idempotents in $A$.
Assume that
$$e_1Ae_2Ae_1=e_1Ae_1\quad\text{and}\quad e_2Ae_1Ae_2=e_2Ae_2.$$
\bnum
\item For any $A$-module $M$, we have
$$e_2Ae_1\otimes_{e_1Ae_1}e_1M\isoto e_2M.$$
\item
$e_1Ae_2$ and $e_2Ae_1$
give a Morita equivalence between $\Mod(e_1Ae_1)$
and $\Mod(e_2Ae_2)$.
\enum
\enlemma
\Proof

\noindent
(i)\ 
The surjectivity follows from
$e_2M=e_2Ae_2M=e_2Ae_1Ae_2M\subset(e_2Ae_1)(e_1M)$.

Let us show its injectivity.
By the assumption, there exists
finitely many elements $a_i\in e_2Ae_1$ and $b_i\in e_1Ae_2$ such that
$e_2=\sum_ia_ib_i$.
Consider now
$u=\sum_jx_j\otimes v_j\in e_2Ae_1\otimes_{e_1Ae_1}e_1M$
(where $x_j\in e_2Ae_1$, $v_j\in e_1M$).
Assume $\sum_jx_jv_j=0$.
Then
$$
u=\sum_{j,i}a_ib_ix_j\otimes v_j
=\sum_{j,i}a_i\otimes b_ix_jv_j=0.
$$
(ii)
It is enough to show that the multiplication maps
$e_2Ae_1\otimes_{e_1Ae_1}e_1Ae_2\to e_2Ae_2$
and $e_1Ae_2\otimes_{e_2Ae_2}e_2Ae_1\to e_1Ae_1$
are isomorphisms. 
For the first one, we apply (i) to $M=Ae_2$. The second one can be handled
similarly.
\QED

The previous result can be expressed in terms of bimodules:
\Prop\label{cor:A}
Let $A$ and $B$ be rings,
and let $P$ be an $(A,B)$-bimodule,
$Q$ a $(B,A)$-bimodule and let
$\varphi\cl P\otimes_B Q\to A$
be a morphism of  $(A,A)$-bimodules,
and $\psi\cl Q\otimes_A P\to B$
a morphism of  $(B,B)$-bimodules.
Assume that $\varphi$ and $\psi$ are surjective and that
the following diagrams commute:
$$\xymatrix@C=8ex{
P\otimes_BQ\otimes_AP\ar[r]^(.6){\varphi\otimes P}\ar[d]^{P\otimes\psi}
&A\otimes_AP\ar[d]^{\can}\\
P\otimes_BB\ar[r]_\can&P}\qquad\text{and}\qquad
\xymatrix@C=8ex{
Q\otimes_AP\otimes_BQ\ar[r]^(.55){\psi\otimes Q}\ar[d]^{Q\otimes \varphi}
&B\otimes_BQ\ar[d]^\can\\
Q\otimes_AA\ar[r]_\can&Q.}
$$
\bnum
\item
Then $\varphi$ and $\psi$ are isomorphisms,
and $P$ and $Q$ give a Morita equivalence between $\Mod(A)$ 
and $\Mod(B)$.
\item
Let $M$ be an $A$-module and $N$ a $B$-module,
and let $f\cl Q\otimes_A M\to N$ and $g\cl P\otimes_BN\to M$
be morphisms such that
the diagrams
$$\xymatrix@C=8ex{
P\otimes_BQ\otimes_AM\ar[r]^(.6){\varphi\otimes M}\ar[d]^{P\otimes f}
&A\otimes_AM\ar[d]^{\can}\\
P\otimes_BN\ar[r]^{g}&M}\qquad\text{and}\qquad
\xymatrix@C=8ex{
Q\otimes_AP\otimes_BN\ar[r]^(.55){\psi\otimes N}\ar[d]^{Q\otimes g}
&B\otimes_BN\ar[d]^\can\\
Q\otimes_AM\ar[r]^{f}&N.}
$$
are commutative. Then $f$ and $g$ are isomorphisms.
\enum
\enprop
\Proof
Apply Lemma \ref{lem:Morg} to
the ring $\left(\begin{matrix} A & P \\ Q & B \end{matrix}\right)$,
its module $\left(\begin{matrix} M\\N\end{matrix}\right)$
and $e_1=\begin{pmatrix}1&0\\0&0\end{pmatrix}$ and
$e_2=\begin{pmatrix}0&0\\0&1\end{pmatrix}$.
\QED

\Remark
\bnum
\item
It would be interesting to describe the image of the morphism
\eqref{eq:Ldelta}.
\item
Let 
$$\Y=\{\frac{m}{d}\ |\ m,d\in\Z,
2\le d\le n, (m,d)=1, m<0\}.$$
It is known that Condition \eqref{eq:HeH} holds for $c\not\in\Y$, while
Condition \eqref{eq:HdiscH} holds when $c-1\not\in\Y$,
cf.\ \cite[Theorem 3.3]{GS}, \cite[Theorem 8.1]{BEG} and
\cite{BE}.
\enum
\enremark

\subsubsection{}

Let us consider the $\D(\ts_\reg)\otimes\D_{\g\times V}$-linear morphism 
\eqn
\sigma\cl\M_c&\to&\M_{c-1}\otimes\det(V) \\
v_c=q(A,z)^c\delta(x,A)&\mapsto& q(A,z)\cdot q(A,z)^{c-1}\delta(x,A)\otimes
l=q(A,z)v_{c-1}\otimes l.
\eneqn
Here $l\in \det(V)\seteq\bigwedge^nV$ is the element such that
$q(A,z)l=A^{n-1}z\wedge A^{n-2}z\wedge\cdots\wedge Az\wedge z$. In
particular, $q(A,z)\otimes l$ 
is a $G$-invariant section of $\CO_{\g\times V}\otimes\det(V)$.

So, the morphism $\sigma$ is $G$-equivariant.
We endow $\M_{c-1}$ with an $H_c$-module structure
via the embedding $\theta_c\cl H_c\hookrightarrow \D(\ts_\reg)\rtimes W$.
Then $\sigma$ is $H_c$-linear.

\Remark
Note that $\M_c\to \M_{c-1}\otimes \det(V)$ is an isomorphism
on $\{q(A,z)\not=0\}$.
However, with our definition of $\M_c$,
the morphism $\M_c\to\M_{c-1}\otimes\det(V)$
is not a monomorphism for certain $c$,
e.g.\ $c=0$. Let us show this after restriction to $\ts_\reg\times V$.
We have $q({}^tA,\partial_z)q(A,z)v_{c-1}=0$ 
for $c=0$ by
\eqref{eq:genM}, while the support of $q({}^tA,\partial_z)v_c$ is the
subvariety $\{q(A,z)=0\}$.
\enremark

Let $\D_{\g\times V}(\disc(x)v_{c-1})$ be the $\D_{\g\times V}$-submodule of
$\M_{c-1}$ generated by $\disc(x)v_{c-1}$.

\Lemma \label{lem:Lcc-1}
\bnum
\item
$\D_{\g\times V}(\disc(x)v_{c-1})$
is invariant by $e_{\det} H_ce_{\det}$.
\item
The morphism $\L_{c-1}\to\D_{\g\times V}(\disc(x)v_{c-1})$
given by $u_{c-1}\mapsto \disc(x)v_{c-1}$ is an isomorphism.
\enum
\enlemma

\sketch
Note that $e_{\det}\disc(x)v_{c-1}=\disc(x)v_{c-1}$.
The proof is similar to that of
Lemma \ref{lem:Lc}: the key point is the following 
(cf.\ e.g.\ \cite[Theorem 3.1]{He})
\eq&&\ysq(\disc(x)v_{c-1})=\Delta_\g(\disc(x)v_{c-1}).\eneq
\QED

By \cite[Proposition 4.1]{BEG},
there is a (unique) isomorphism
$$f\cl e_{\det}H_ce_{\det}\iso eH_{c-1}e$$
such that $\theta_{c-1}\bl(f(a)\br)=\disc(x)^{-1}\theta_c(a)\disc(x)$ for
$a\in e_{\det}H_ce_{\det}$.

The isomorphism $\L_{c-1}\iso\D_{\g\times V}(\disc(x)v_{c-1})$
of Lemma \ref{lem:Lcc-1} is compatible with $f$ 
and we will sometimes view
$\L_{c-1}$ as an
$(e_{\det}H_ce_{\det}\otimes\D_{\g\times V})$-module.

\medskip
By Lemma \ref{lem:reg}, the image of the morphism
$$e_{\det}H_ce\otimes_{eH_ce} \L_c\vert_{\g_\reg\times V}\to
\M_c\vert_{\g_\reg\times V},\ a\otimes u_c\mapsto av_c$$
is contained in $\D_{\g_\reg\times V}(\disc(x)v_c)$.
It follows from Lemma \ref{lem:Lcc-1} that over $\g_\reg\times V$,
the composite morphism
$e_{\det}H_ce\otimes_{eH_ce} \L_c\to\M_c
\to\M_{c-1}\otimes \det(V)$
factors through a morphism
\eq
&&\vphi\cl e_{\det}H_ce\otimes_{eH_ce} \L_c\vert_{\g_\reg\times V}
\To \L_{c-1}\otimes \det(V)\vert_{\g_\reg\times V}.
\label{arr:detp}
\eneq

Similarly, we have
the morphism 
\begin{equation}
\label{mor:Lc-1c}
\begin{split}
\psi\cl eH_ce_{\det}\otimes_{e_{\det}H_ce_{\det}}
\L_{c-1}\otimes\det(V)\vert_{\{q(A,z)\not=0\}}
&\to \L_{c}\vert_{\{q(A,z)\not=0\}}\\
a\otimes u_{c-1}\otimes l&\mapsto (a\disc(x))q(A,z)^{-1}u_c.
\end{split}
\end{equation}

The morphism $\varphi$ is linear over
$e_{\det}H_ce_{\det}\simeq eH_{c-1}e$ and
the morphism $\psi$ is linear over $eH_ce$.
We have
$$\vphi(\disc(x)e\otimes u_c)=q(A,z)u_{c-1}\otimes l$$
and
$$q(A,z)\psi(\disc(x)e_{\det}\otimes u_{c-1}\otimes l)=
\disc^2(A)u_c$$
where $\disc^2(A)$ is the discriminant of the characteristic polynomial of $A$.

Note that the following diagrams commute on 
$\g_\reg\times V\cap\{q(A,z)\not=0\}$:
\eq
&&\ba{c}
\xymatrix@C=8ex{
eH_ce_{\det}\tens_{e_{\det}H_ce_{\det}}e_{\det}H_ce\tens_{eH_ce}\L_c
\ar@<.6ex>[r]^(.6){}\ar[d]^{\varphi}
&eH_ce\tens_{eH_ce}\L_c\ar[d]^{\can}\\
eH_ce_{\det}\tens_{e_{\det}H_ce_{\det}}\bl(\L_{c-1}\otimes\det(V)\br)
\ar@<.5ex>[r]^(.7){\psi}&\L_c}
\ea\label{dia:Lc}
\eneq
and
\eq&&\ba{c}\xymatrix@C=3ex{
e_{\det}H_ce\tens_{eH_ce}eH_ce_{\det}\tens_{e_{\det}H_ce_{\det}}
\bl(\L_{c-1}\otimes\det(V)\br)\ar@<.6ex>[r]^(.55){}\ar[d]^{\psi}
&e_{\det}H_ce_{\det}\tens_{e_{\det}H_ce_{\det}}\bl(\L_{c-1}\otimes\det(V)\br)
\ar[d]^\can\\
e_{\det}H_ce\tens_{eH_ce}\L_c\ar@<.5ex>[r]^{\varphi}&\L_{c-1}\otimes\det(V).}
\ea\label{dia:Lc-1}
\eneq

\Prop\label{prop:phiext}
The morphism
$\vphi$ extends uniquely to
a morphism of $\D_{\g\times V}$-modules:
\eq
&&\vphi\cl e_{\det}H_ce\otimes_{eH_ce} \L_c
\To \L_{c-1}\otimes \det(V).
\label{arr:det}
\eneq
\enprop

The proof will proceed by reduction to rank two. Recall that 
$\g_1$ denotes the open subset of $\g$
of matrices with at least $(n-1)$ distinct eigenvalues.
Then $\g\setminus\g_1$ is a closed subset of $\g$ of codimension $2$.

We shall prove first the following lemma.
\Lemma
\label{redrank2}
After restriction to $\g_1\times V$, we have an inclusion of submodules
of $\M_{c-1}$
$$H_c\D_{\g\times V}\ol{v}_c
\subset \C[\ts]\D_{\g\times V}\ol{v}_c+\C[\ts]\D_{\g\times V}\disc(x)v_{c-1}$$
where $\ol{v}_c=q(A,z)v_{c-1}$.
\enlemma
\sketch
Since $H_c=\C[\ts]\C[\ts^*]\C[W]$,
it is enough to show that
\eq
&&\C[\ts^*]\D_{\g\times V}\ol{v}_c\in  
\C[\ts]\D_{\g\times V}\ol{v}_c+\C[\ts]\D_{\g\times V}\disc(x)v_{c-1}
\quad\text{on $\g_1\times V$.}
\label{eq:ctLc}
\eneq
Here the action of $\C[\ts^*]$ is through
$\C[\ts^*]\into H_c\To[{\theta_c}]\D(\ts_\reg)\rtimes W$.

\smallskip
Let us assume first that $n=2$.
We have
$$q(A,z)=-A_{21}z_1^2+(A_{11}-A_{22})z_1z_2+A_{12}z_2^2.$$
We put
$$q(\partial_A,z)=-z_1^2\partial_{A_{12}}+z_1z_2(\partial_{A_{11}}-
\partial_{A_{22}})+ z_2^2\partial_{A_{21}}.$$
We will show that 
\eq
(\partial_{x_1}-\partial_{x_2})q(A,z)v_{c-1}
=-q(\partial_A,z)(x_1-x_2) v_{c-1}.
\label{shiftrank2}
\eneq
This is an equality in the $\D_{\g\times V}$-submodule $\iota(\L_{c-1})$
of $\M_{c-1}$. 
Note that $(y_1-y_2)v_{c-1}=(\partial_{x_1}-\partial_{x_2})v_{c-1}$.

By \S\,\ref{subsec:delta},
we have $$v_{c-1}=q(A,z)^{c-1}\delta(x_1+x_2-\tr(A))\delta(x_1x_2-\det(A)).$$
Since $q(\partial_A,z)q(A,z)=q(\partial_A,z)\tr(A)=0$ and
$q(\partial_A,z)\det(A)=-q(A,z)$, we obtain
$$q(\partial_A,z)v_{c-1}
=q(A,z)^{c}\delta(x_1+x_2-\tr(A))\delta'(x_1x_2-\det(A)).$$
On the other hand, we have
$$(\partial_{x_1}-\partial_{x_2})q(A,z)v_{c-1}
=(x_2-x_1)q(A,z)^{c}\delta(x_1+x_2-\tr(A))\delta'(x_1x_2-\det(A)).$$
The equality \eqref{shiftrank2} then follows.

\medskip
We assume now $n\ge 2$.
Let $S$ be the locally closed subset of $\g$ of matrices
$$\begin{pmatrix}A'&0&0&\cdots\\0&a_3&0&\cdots\\
0&0&a_4\\
\bxes{-.65ex}
&\bxes{-.65ex}&&
\bxes{.3ex}
\\
&&&&a_n
\end{pmatrix}$$
where $A'$ is a $2\times 2$ matrix
and $a_i\not=a_j$ ($3\le i<j\le n$) and $a_i$
is not an eigenvalue of $A'$ for $3\le i\le n$.
Let $\ts_1=\ts\cap S=\set{x\in\ts}{\text{$x_i\not=x_j$ for
$i<j$ and $3\le j$}}.$
Let $x'=(x_1,x_2)$, $x''=(x_3,\ldots,x_n)$ and $a''=(a_3,\ldots,a_n)$.

We have $G\cdot S=\g_1$.
Let $i\cl S\times V\into\g\times V$ be 
the inclusion map. Then,
$i$
is non-characteristic for $\L_c$ and $\M_{c-1}$,
because we have $T_xS+T_x(G\cdot x)=T_x\g$ for any $x\in S$.

\smallskip
Denote by $\g'$ the subalgebra of $\g$
of matrices $(A_{ij})$ with $A_{ij}=0$ whenever $i>2$ or $j>2$.
We identify $\g'$ with $\gl_2(\C)$. Given an object $\mathcal{X}$ defined
earlier for $\g$, we denote by $\mathcal{X}'$ the corresponding
objects for $\g'$ (i.e., case $n=2$). For example, $W'$ is the subgroup of
$W$ generated by $s_{12}$.

\smallskip
Let $i''\cl\ts\times S\to\ts\times\g$ be the embedding.
We have an isomorphism of $\D_{\ts\times S}$-modules compatible with
the action of $W$ (cf.\ Proof of Lemma \ref{lem:reg}):
$$i^{\prime\prime *}\D_{\ts\times\g}\delta(x,A)
\xrightarrow[\sim]{\delta(x,A)\mapsto\sum_w  T_w^* \disc_1(A',a'')^{-1}
\delta(x',A')\delta(x''-a'')}
\bigoplus_{w\in W'\backslash W} T_w^* \D_{\ts\times S} \delta(x',A')
\delta(x''-a'').$$
Here, $T_w$ is the automorphism of $\ts$ given by $w$,
and $\disc_1(A',a'')=\disc(a'')\prod_{i=3}^n\det(a_iI_2-A')$,
$\delta(x',A')=\delta(x_1+x_2-\tr(A'))\delta(x_1x_2-\det(A'))$.

\smallskip
Let $A\in S$.
We have
$$q(A,z)=q'(A',z')\cdot q_1(A,z),$$
where
$$q_1(A,z)=(z_3\cdots z_n)\disc_1(A',a'').$$
Note that $\disc_1(A',a'')$ is invertible on $S$.

Let $p\cl\ts_\reg\times S\times V\to S\times V$ be the projection.
We have a $\D(\ts_\reg)\otimes \D_{S\times V}$-linear 
isomorphism compatible with
the action of $W$:
\eq
\label{eqisorank2}
i^*\M_c\xrightarrow[\sim]{v_c\mapsto e\otimes \widetilde{v}_c}
\C[W]\otimes_{\C[W']} p_*\bigl(\D_{\ts_\reg\times S\times V}\widetilde{v}_c\bigr)
\eneq
where $\widetilde{v}_c=v'_c\, q_1(A,z)^c\,\disc_1(A',a'')^{-1}\,
\delta(x''-a'')$ with
$v'_c=q'(A',z')^{c}\,\delta(x',A')$.
Note that $s_{12}$ acts trivially on $\widetilde{v}_c$.
The action of $\D(\ts_\reg)\rtimes W$ on $\C[W]\otimes_{\C[W']}
p_*\bigl(\D_{\ts_\reg\times S\times V}\widetilde{v}_c\bigr)$
is given by:
$$(a\otimes w)(w'\otimes s)=
(ww')\otimes \bigl(((ww')^{-1}a)s\bigr)\quad
\text{for $w$, $w'\in W$, $a\in \D(\ts_\reg)$, 
$s\in p_*\bigl(\D_{\ts_\reg\times S\times V}\widetilde{v}_c\bigr)$.}$$
Note that
$\D_{S\times V}\widetilde{v}_c$ is stable by $\C[\ts_1]^{W'}$
as a submodule of 
$p_*\bigl(\D_{\ts_\reg\times S\times V}\widetilde{v}_c\bigr)$.
Since $\C[\ts_1]=\C[\ts]\C[\ts_1]^{W'}$,
$\C[\ts]\D_{S\times V}\widetilde{v}_c$ is stable by $\C[\ts_1]$.

\smallskip
Let us still denote by $\widetilde{v}_c=q(A,z)\widetilde{v}_{c-1}$,
the image of $\widetilde{v}_c$.

Let us set $\widetilde{y}_1=\partial_{x_1}-c(x_1-x_2)^{-1}(1-s_{12})$
and $\widetilde{y}_2=\partial_{x_2}-c(x_2-x_1)^{-1}(1-s_{12})$,
partial Dunkl operators,
and let $R$ be the algebra generated by $\widetilde{y}_1$,
$\widetilde{y}_2$ and $\partial_{x_i}$ ($i=3,\ldots,n$).
Then $s_{12}$ acts on $R$ by the permutation of $\widetilde{y}_1$
and $\widetilde{y}_2$.
We have $R=R^{W'}\oplus (\widetilde{y}_1-\widetilde{y}_2)R^{W'}$.

Let
\eqn\tN
&=&\C[\ts]\D_{S\times V}\widetilde{v}_c+
(\widetilde{y}_{1}-\widetilde{y}_{2})\C[\ts]\D_{S\times V}\widetilde{v}_c\\
&=&\C[\ts]\D_{S\times V}\widetilde{v}_c+
\C[\ts]\D_{S\times V}(\widetilde{y}_{1}-\widetilde{y}_{2})\widetilde{v}_c\\
&=&\C[\ts]\D_{S\times V}\widetilde{v}_c
+\C[\ts]\D_{S\times V}(\partial_{x_1}-\partial_{x_2})\widetilde{v}_{c}
\eneqn
be a submodule of
$p_*\Bigl(\D_{\ts_\reg\times S\times V}\widetilde{v}_{c-1}\Bigr)$.
Since $(\widetilde{y}_1+\widetilde{y}_2)\widetilde{v}_c$,
$\widetilde{y}_1\widetilde{y}_2\widetilde{v}_c$, and $\partial_{x_i}\widetilde{v}_c$
($i=3,\ldots,n$)
belong to $\C[\ts]\D_{S\times V}\widetilde{v}_c$ (cf.\ Lemma~\ref{lem:Lc}),
$\tN$ is invariant by $R$.

\smallskip
Set
$\N=\C[W]\otimes_{\C[W']}\tN$.
Let us show that $\N$ is invariant by the action of
$\C[\ts^*]\subset H_c\subset\D(\ts_\reg)\rtimes W$.
For any $i$, we have 
$$y_i(w\otimes t)
=w\otimes \partial_{x_{w^{-1}(i)}}t-
c\sum_{k\not=i}w(1+s_{{w^{-1}(i)},{w^{-1}(k)}})
\otimes (x_{w^{-1}(i)}-x_{w^{-1}(k)})^{-1}t$$
for any $w\in W$ and $t\in \tN$.
Since $(x_a-x_b)^{-1}\in\C[\ts_1]$ when $a$ or $b$
is in $\{3,\ldots,n\}$,
we have
$y_i(w\otimes t)\in\N$ when $w^{-1}(i)\not=1,2$.
If $w^{-1}(i)=1$, then
\eqn
y_i(w\otimes t)
&\equiv& w\otimes \partial_{x_1}t
-cw(1+s_{12})\otimes(x_1-x_2)^{-1}t\ \mod \N\\
&=&w\otimes \widetilde{y}_1t\in \N.
\eneqn
The case $w^{-1}(i)=2$ is similar.
Hence we have shown that $\N$ is invariant by $\C[\ts^*]$.
Thus, we obtain
\eqn
\C[\ts^*](e\otimes\widetilde{v}_c)\subset\N.
\eneqn

\smallskip
The study of rank $2$ above, i.e.\ \eqref{shiftrank2}, shows that 
$$(\widetilde{y}_1-\widetilde{y}_2)\widetilde{v}_c
\subset \C[\ts]\D_{S\times V}\widetilde{v}_c
+\C[\ts]\D_{S\times V}(x_1-x_2)\widetilde{v}_{c-1}.$$
Hence we obtain
$$\tN\subset\tN'\seteq\C[\ts]\D_{S\times V}\widetilde{v}_c
+\C[\ts]\D_{S\times V}\disc(x)\widetilde{v}_{c-1},$$
which implies
\eq
\C[\ts^*](e\otimes\widetilde{v}_c)
\subset \N'\seteq\C[W]\otimes\tN'.
\label{eq:bvN}
\eneq

We have a commutative diagram, where the horizontal map is an
isomorphism
$$\xymatrix{
W\times_{W'}\ts_1 \ar[rrrr]^{(w,x)\mapsto (w(x),x)}_\sim
\ar[drr]_{(w,x)\mapsto w(x)} &&&&
\ts\times_{\ts/W}\ts_1/W' \ar[dll]^{(x,x')\mapsto x} \\
&& \ts}$$
The diagram above is $W$-equivariant, for the action of $g\in W$ given by
$$g\cdot(w,x)=(gw,x) \text{ for }(w,x)\in W\times_{W'}\ts_1$$
$$g\cdot(x,x')=(g(x),x') \text{ for } (x,x')\in \ts\times_{\ts/W}\ts_1/W'.$$
It follows that we have an isomorphism of $\C[\ts]$-modules
\begin{align*}
\C[\ts]\otimes_{\C[\ts]^W} \C[\ts_1]^{W'}&\iso
\C[W]\otimes_{\C[W']}\C[\ts_1] \\
a\otimes a'&\mapsto \sum_{w\in W/W'} w\otimes w^{-1}(a)a'.
\end{align*}
In particular, we have
$\C[W]\otimes_{\C[W']}\C[\ts_1]=\C[\ts]\cdot(e\otimes\C[\ts_1]^{W'})$.
Since $\C[\ts_1]^{W'}\widetilde{v}_c\subset \D_{S\times V}\widetilde{v}_c$
and $\C[\ts_1]^{W'}\disc(x)\widetilde{v}_{c-1}\subset 
\D_{S\times V}\disc(x)\widetilde{v}_{c-1}$,
we deduce that
$$\N'=\C[\ts]\Bigl(
e\otimes\D_{S\times V}\widetilde{v}_c
+e\otimes\D_{S\times V}\disc(x)\widetilde{v}_{c-1}\Bigr).$$
Together with \eqref{eq:bvN}, we obtain
$$\C[\ts^*]\D_{S\times V}(e\otimes\widetilde{v}_c)
\subset \C[\ts]\Bigl(
\D_{S\times V}(e\otimes\widetilde{v}_c)
+\D_{S\times V}(e\otimes\disc(x)\widetilde{v}_{c-1})\Bigr).$$
Via the isomorphism
\eqref{eqisorank2}, this shows that
$$i^*\Bigl(\C[\ts^*]\D_{\g\times V}\ol{v}_c\Bigr)
\subset
i^*\Bigl(\C[\ts]\D_{\g\times V}\ol{v}_c+\C[\ts]\D_{\g\times V}\disc(x)v_{c-1}\Bigr).$$
Since $\mu^{-1}(0)\cap T^*_{S\times V}(\g\times V)\subset
T^*_{\g\times V}(\g\times V)$, the non-characteristic condition implies
the desired result \eqref{eq:ctLc} (cf.\ \S\,\ref{subsubnonchar}).
\QED

\begin{proof}[Proof of Proposition \ref{prop:phiext}]
By Lemma \ref{redrank2}, we have, on $\g_1\times V$,
\eqn
e_{\det}H_c\D_{\g\times V}\ol{v}_c
&\subset& e_{\det}\C[\ts]\D_{\g\times V}\ol{v}_c+
e_{\det}\C[\ts]\D_{\g\times V}\disc(x)v_{c-1}\\
&\subset&\C[\ts]^W\disc(x)\D_{\g\times V}\ol{v}_c+
\C[\ts]^W\D_{\g\times V}\disc(x)v_{c-1}
=\D_{\g\times V}\disc(x)v_{c-1},
\eneqn
since $e_{\det}\C[\ts]e=\C[\ts]^W\disc(x)e$ and $e_{\det}\C[\ts]e_{\det}=
\C[\ts]^We_{\det}$.
Hence $\vphi$ extends to a morphism defined
on $\g_1\times V$.
Then the desired result follows from
$\SH^1_{(\g\setminus\g_1)\times V}(\L_{c-1})=0$ (Lemma \ref{vanishcoh}).
\end{proof}

\section{Cherednik algebras and Hilbert schemes}
\label{secHilb}
\subsection{Geometry of the Hilbert scheme}
\label{secgeo}
\subsubsection{}
We refer to \cite{N, Ha1} for basic results on Hilbert schemes of points on
$\C^2$.

Let us recall that
$$\X=
\set{(A,B,z,\zeta)\in\g\times\g\times V\times V^*}{\C\langle A,B\rangle z=V}$$
is the set of stable points for the action of $G$ on
$T^*(\g\times V)$, relative to the
character $\det$ of $G$.
The group $G$ acts freely on $\X$.
Let $\mu_\X\cl \X\to\g$ be the moment map:
$$\mu_\X(A,B,z,\zeta)=-[A,B]-z\circ\zeta.$$
It is a smooth morphism.
Let $\Hilb^n(\C^2)$ be the Hilbert scheme
classifying closed subschemes of $\C^2$ with length $n$.
Then we have an isomorphism
$\Hilb^n(\C^2)\iso \mu_\X^{-1}(0)/G$.
Note that we have $\zeta=0$ on $\mu_\X^{-1}(0)$
(cf.\ \cite[Lemma 2.3]{GG}).

We shall write $\Hilb$ instead of $\Hilb^n(\C^2)$ for short.
Let us denote by
$p\cl \mu_\X^{-1}(0)\to \Hilb$
the quotient map.

Let us recall the construction of $p$.
For $(A,B,z,\zeta)\in\mu_\X^{-1}(0)$,
we regard $V$ as a $\C[X,Y]$-module
by $X\mapsto A$ and $Y\mapsto B$.
Then $z$ gives an epimorphism 
$\C[X,Y]\epi V$ of $\C[X,Y]$-modules.
Hence $V$ gives a closed subscheme of $\C^2=\Spec(\C[X,Y])$
of length $n$, which is the corresponding point
of $\Hilb$.

\smallskip
Let $\pi\cl\Hilb\to(\ts\times\ts^*)/W$
be the Hilbert-Chow morphism.
Then $\Hilb$ is a resolution of singularities of $(\ts\times\ts^*)/W
\simeq (\C^2)^n/S_n$,
the scheme of
$n$ unordered points in $\C^2$.
We have canonical isomorphisms
$$\Gamma(\mu_{\X}^{-1}(0),\CO_{\mu_{\X}^{-1}(0)})^G\iso\Gamma(\Hilb,\CO_\Hilb)
\iso \Gamma((\ts\times\ts^*)/W, \CO_{(\ts\times\ts^*)/W})
\iso\C[\ts\times\ts^*]^W.$$
Let $(\ts\times\ts^*)_\reg$ be the open subset of $\ts\times\ts^*$ 
where the action of $W$ is free.
The Hilbert-Chow morphism $\pi$ is an isomorphism over
$(\ts\times\ts^*)_\reg/W$.
Let $E\seteq\pi^{-1}\Bigl(\bl((\ts\times\ts^*)
\setminus(\ts\times\ts^*)_\reg\br)/W\Bigr)$
be the exceptional divisor.
It is a closed irreducible hypersurface of $\Hilb$.
The line bundle $L$ on $\Hilb$
associated with the $G$-equivariant line bundle 
$\CO_\X\otimes\det(V)$ on $\X$
is a very ample line bundle on $\Hilb$.

Let us set
$$\C[\mu_\X^{-1}(0)]^{G,\det}
=\set{\phi(p)\in\C[\mu_\X^{-1}(0)]}
{\text{$\phi(gp)=\det(g)\phi(p)$ for any $g\in G$}}.$$
It is isomorphic to $\Gamma(\Hilb,L)
\simeq(\C[\mu_\X^{-1}(0)]\otimes\det(V))^G$.
Let $i\cl \ts\times\ts^*\times V\hookrightarrow \g\times\g\times V\times V^*$
be the embedding with the last component $\zeta=0$.
Then $i^{-1}(\mu_\X^{-1}(0))$
contains $(\ts_\reg\times \ts^*
\cup \ts\times \ts_\reg^*)\times (\C^*)^n$.
For any $\phi\in \C[\mu_\X^{-1}(0)]^{G,\det}$,
we have
$(i^*\phi)(x,y,gz)=\det(g)(i^*\phi)(x,y,z)$
for any invertible diagonal matrix $g$.
Hence we have
$$(i^*\phi)(x,y,z)=a(x,y)(z_1\cdots z_n)$$
for some rational function $a(x,y)$ 
which is regular on
$(\ts_\reg\times \ts^*)\cup(\ts\times \ts_\reg^*)$, an open
subset of $\ts\times\ts^*$ with complement of codimension $2$.
Hence we have
$$a(x,y)\in\C[\ts\times\ts^*]^{W,\det}=
\set{a\in \C[\ts\times\ts^*]}{\text{$wa=\det(w)a$ for any $w\in W$}}.$$
Thus we obtain a map which is known to be an
isomorphism (cf.\ e.g. \cite[Proposition 8.2.1]{GG})
and
we denote its inverse by $i_d$:
\begin{equation}
\label{mor:id}
\begin{split}
\C[\mu_\X^{-1}(0)]^{G,\det}\otimes\det(V)&\iso \C[\ts\times\ts^*]^{W,\det}\\
\phi\otimes l&\mapsto\langle l,z_1\wedge\cdots\wedge z_n\rangle a.
\end{split}
\end{equation}

Similarly, we have an isomorphism (cf.\ e.g. \cite[Lemma 2.7.3]{GG})
whose inverse we denote by $i_s$:
\eq
\C[\mu_\X^{-1}(0)]^G\iso\C[\ts\times\ts^*]^{W}.
\label{mor:is}\eneq

Summarizing, we have the following isomorphisms
\eq
&&\ba{rcccccl}
i_d&\cl& \C[\ts\times\ts^*]^{W,\det}&\isoto&\C[\mu_\X^{-1}(0)]^{G,\det}
\otimes\det(V)&\simeq&\Gamma(\Hilb,L),\\[2ex]
i_s&\cl &\C[\ts\times\ts^*]^{W}&\isoto&\C[\mu_\X^{-1}(0)]^G
&\simeq&\CO_\Hilb(\Hilb).
\ea
\eneq

\subsubsection{}
For a subset $Y$ of $\Z_{\ge0}\times\Z_{\ge0}$ with cardinality $n$,
set $p_Y=\det(x_k^iy_k^j)_{(i,j)\in Y,\,k=1,\ldots n}
\in\C[\ts\times\ts^*]^{W,\det}$ and
$s_Y(A,B,z,\zeta)=\det(A^iB^jz)_{(i,j)\in Y}\in\C[\mu_\X^{-1}(0)]^{G,\det}
=L(\Hilb)$.
Then $\{p_Y\}_Y$ is a basis of $\C[\ts\times\ts^*]^{W,\det}$
as a vector space and
$i_d(p_Y)=s_Y$.
The $\CO_\Hilb$-module $L$ is generated by 
$\{s_Y\}_{Y}$,
where $Y$ ranges over the set of Young diagrams of size $n$.
Here we regard a Young diagram $Y$ as a subset
of $\Z_{\ge0}\times \Z_{\ge0}$
such that $(i,j)\in Y$ as soon as $(i,j+1)$ or $(i+1,j)$
belongs to $Y$.

There is a canonical global section $\tau\in
\Gamma(\Hilb;L^{\otimes -2})$
satisfying the following property:
\eq
i_d(a_1)i_d(a_2)\tau=i_s(a_1a_2)
\quad\text{for any $a_1$, $a_2\in \C[\ts\times\ts^*]^{W,\det}$.}
\eneq
Note that $\tau$ is identified with a function on $\mu_\X^{-1}(0)$
such that $\tau(gp)=\det(g)^{-2}\tau(p)$ ($p\in \mu_\X^{-1}(0)$
and $g\in G$).

The exceptional divisor $E$ coincides with the set of zeroes of $\tau$,
and we obtain an isomorphism
$$L^{\otimes 2}\iso \CO_\Hilb(-E).$$
Let us denote by $\disc^2(A)$ the discriminant of the characteristic
polynomial of $A$,
and similarly for $\disc^2(B)$.
Then we have
\eqn
i_d(\disc(x))=q(A,z),\quad i_d(\disc(y))=q(B,z),
\quad i_s(\disc(x)^2)=\disc^2(A),\quad i_s(\disc(y)^2)=\disc^2(B).
\eneqn
Hence we have
$$\text{$\disc^2(A)=q(A,z)^2\tau$ and
$\disc^2(B)=q(B,z)^2\tau$.}
$$

\Lemma
\label{Hilbcodim2}
\bnum
\item
The hypersurface of $\mu_\X^{-1}(0)$ defined by 
$q(A,z)=0$ is irreducible,
and $p^{-1}E\cap\{q(A,z)=0\}$ is of codimension $2$ in $\mu_\X^{-1}(0)$.
\item
The hypersurface of $\mu_\X^{-1}(0)$ defined by $\disc^2(A)=0$ is
$p^{-1}E\cup\{q(A,z)=0\}$.
\item
$\mu_\X^{-1}(0)\cap
\{q(A,z)=q(B,z)=0\}$ is of codimension $2$
in $\mu_\X^{-1}(0)$.
\enum
\enlemma
Note that (i) follows from the fact that
$q(A,z)$ does not vanish on the irreducible hypersurface $p^{-1}E$
of $\mu_\X^{-1}(0)$, and $q(A,z)$ is irreducible on 
$\mu_\X^{-1}(0)\setminus p^{-1}E$. Statement (iii) follows from
\cite[Lemma 3.6.2]{Ha1}.

\subsection{W-algebras on the Hilbert scheme}
\label{subsec:Ana}
\subsubsection{}
In the preceding sections, we have regarded $\X$, $\Hilb$, etc.\ as schemes.
Hereafter, we regard them as complex manifolds. Note that the 
previous constructions
and results would remain valid in the analytic category.
Let $\W_\X$ be the $\W$-algebra on $\X$ associated with
$\D_{\g\times V}$.
Denoting by $\pi\cl \X\to \g\times V$ the projection,
we have a ring homomorphism
$\pi^{-1}\D_{\g\times V}\to\W_\X$
respecting the order filtration.
The ring $\W_\X$ is flat over $\pi^{-1}\D_{\g\times V}$.
The action of $G$ on $\g\times V$ induces an 
action of $G$ on $\W_\X$ and there is a 
quantized moment map
$\mu_\W\cl\g\to\W_\X$.

We have morphisms
$$\kappa_0\cl\C[\ts]^W\isoto\C[\g]^G\to\W_\X(\X)$$
and
$$\kappa_1\cl\C[\Gt^*]^W\isoto\C[\g^*]^G\hookrightarrow
\D_\g(\g)\to\W_\X(\X).$$
Note that $\kappa_1(\ysq)=\Delta_\g$.

For $k\in\Z_{\ge0}$, let $\C[\Gt^*]^W_k$ be the homogeneous part of
$\C[\Gt^*]^W$ of degree $k$.
Then $\kappa_0$ sends $\C[\ts]^W$ to $\W_\X(0)$ and
$\kappa_1$ sends $\C[\Gt^*]^W_k$ to $\W_\X(k)$ and we have 
the commutative diagrams:
\eq
\ba{c}
\xymatrix{
\C[\ts]^W\ar[r]^{\kappa_0}\ar[dr]&\W_\X(0)\ar[d]^{\sigma_0}\\
&{\CO_\X}}\ea\quad\text{and}\quad
\ba{c}
\xymatrix{
\C[\Gt^*]^W_k\ar[r]^{\kappa_1}\ar[dr]_{\h^{-k}}&\W_\X(k)\ar[d]^{\sigma_k}\\
&{\h^{-k}\CO_\X}}
\ea\label{eq:symb}
\eneq
Let us consider
$\W_\X\otimes_{\D_{\g\times V}}\L_c$,
which we denote by the same letter $\L_c$.
With the notation of \S\,\ref{sectwisted},
we have $\L_c=\Phi_{c\tr}(\W_\X)$.
Hence $\L_c$ is a twisted $G$-equivariant $\W_\X$-module with twist
$c\tr$. 
Let $u_c$ be the canonical section of $\L_c$,
and set $\L_c(m)=\W_\X(m)u_c$.
Then we have an isomorphism
$$\L_c(0)/\L_c(-1)\iso\CO_{\mu_\X^{-1}(0)}.$$
The support of $\L_c$ is $\mu_\X^{-1}(0)$.
The $\W_\X$-module $\L_c$ has a left action of $eH_ce$ by Lemma \ref{lem:Lc}.
Via the anti-involution $h\mapsto h^*$ of $H_c$,
we regard $\L_c$ as a $(\W_\X,eH_ce)$-bimodule.
Similarly,
$\L_{c-1}$ has a structure of $(\W_\X,e_{\det}H_ce_{\det})$-bimodule
(Lemma \ref{lem:Lcc-1}).
These actions are explicitly given by:
\eq
&&\ba{l}
u_cea=\kappa_0(a)u_c \quad\text{for $a\in\C[\ts]^W\subset H_c$,}\\[1ex]
u_ceb=\kappa_1(b)u_c \quad\text{for $b\in\C[\Gt^*]^W\subset H_c$.}
\ea\\[2ex]
&&\ba{l}
u_{c-1}e_{\det}a=\kappa_0(a)u_{c -1}
\quad\text{for $a\in\C[\ts]^W\subset H_c$,}\\[1ex]
u_{c-1}e_{\det}b=\kappa_1(b)u_{c-1}
 \quad\text{for $b\in\C[\Gt^*]^W\subset H_c$.}
\ea
\eneq

\smallskip
Since $\mu_\X^{-1}(0)$ is smooth, we have
$$\text{$\ext_{\W_\X}^j(\L_c,\W_\X)=0$ for $j\not=\codim_\X(\mu_\X^{-1}(0))$.}
$$
Hence,
for any closed subset $S\subset\mu_\X^{-1}(0)$,
we have by Lemma~\ref{lem:vanloc}:
\eq
\text{$\SH^j_S(\L_c)=0$ for $j<\codim_{\mu_\X^{-1}(0)}S$.}
\label{eq:ext}
\eneq

In \eqref{arr:det} and \eqref{mor:Lc-1c},
we defined the morphisms:
\eq
&&\vphi\cl \L_c \tensor_{eH_ce}
eH_ce_{\det}
\To \L_{c-1}\otimes \det(V)
\eneq
and
$$\psi\cl \bl(\L_{c-1}\otimes \det(V)\br) \tensor_{e_{\det}H_ce_{\det}}
e_{\det}H_ce\mid_{\{q(A,z)\not=0\}}
\To \L_c\mid_{\{q(A,z)\not=0\}}.
$$

\Prop\label{prop:psiext}
The morphism $\psi$ extends uniquely to a morphism defined on $\X$.
\enprop
\sketch
We have 
$$q(A,z)\psi(u_{c-1}\otimes a)=u_c\cdot(\disc(x)a)$$
for any $a\in e_{\det}H_ce$.

Now let us show that
\eq
&&(\ad(\Delta_\g)^kq(A,z))\psi(u_{c-1}\otimes a)
=u_c\cdot((\ad(\ysq)^k\disc(x))a)
\label{eq:ucuc-1}
\eneq
holds on $\{q(A,z)\ne0\}$ by the induction on $k$.

We have
\eqn
&&(\ad(\Delta_\g)^kq(A,z))\psi(u_{c-1}\otimes a)\\
&&\hs{5ex}=\Delta_\g (\ad(\Delta_\g)^{k-1}q(A,z))\psi(u_{c-1}\otimes a)
-(\ad(\Delta_\g)^{k-1}q(A,z))\Delta_\g \psi(u_{c-1}\otimes a).
\eneqn
The first term is calculated as
\eqn
\Delta_\g (\ad(\Delta_\g)^{k-1}q(A,z))\psi(u_{c-1}\otimes a)
&=&\Delta_\g u_c\cdot((\ad(\ysq)^{k-1}\disc(x))a)\\
&=&u_c\ysq\cdot((\ad(\ysq)^{k-1}\disc(x))a)\\
&=&u_c\cdot(\ysq(\ad(\ysq)^{k-1}\disc(x))a).
\eneqn
The second term is calculated as
\eqn
(\ad(\Delta_\g)^{k-1}q(A,z))\Delta_\g \psi(u_{c-1}\otimes a)
&=&(\ad(\Delta_\g)^{k-1}q(A,z))\psi(\Delta_\g u_{c-1}\otimes a)\\
&=&(\ad(\Delta_\g)^{k-1}q(A,z))\psi(u_{c-1}\ysq \otimes a)\\
&=&(\ad(\Delta_\g)^{k-1}q(A,z))\psi(u_{c-1}\otimes \ysq a)\\
&=&u_{c}\cdot((\ad(\ysq)^{k-1}\disc(x))\ysq a).
\eneqn
Hence we obtain \eqref{eq:ucuc-1}.
In particular, letting $k$ to be 
$n(n-1)/2$, the degree of $\disc(x)$, and using the fact that
$\ad(\Delta_\g)^{n(n-1)/2}q(A,z)$ is equal to $q(\partial_A,z)$ 
up to a constant multiple (see e.g.\ \eqref{eq:ad} and the sentence below),
we obtain
\eq
&&q(\partial_A,z)\psi(u_{c-1}\otimes a)=
u_c\cdot(\disc(y)a).
\eneq
Hence $\psi(u_{c-1}\otimes a)$ extends to a section of
$\L_c$ outside $q(B,z)=0$.

Thus we have shown that
$\psi(u_{c-1}\otimes a)$ is a section defined outside
$\{q(A,z)=0\}\cap \{q(B,z)=0\}$.
Since $\{q(A,z)=0\}\cap \{q(B,z)=0\}\cap \mu_\X^{-1}(0)$
is of codimension $2$ in $\mu_\X^{-1}(0)$ (Lemma \ref{Hilbcodim2}),
it follows that
$\psi(u_{c-1}\otimes a)$ extends to a global section of
$\L_c$ by \eqref{eq:ext}.
\QED

\Remark\label{rem:4.3}
\bnum
\item
So, we have obtained a structure of
$(e+e_{\det})H_c(e+e_{\det})$-module on
$\L_c\oplus\L_{c-1}\otimes\det(V)$.
\item
We have
\eq
&&\ba{l}
\ba{l}
\vphi(u_c\otimes e\disc(x))=q(A,z)u_{c-1},\\[1ex]
\vphi(u_c\otimes e\disc(y))=q(\partial_A,z)u_{c-1},
\ea\\[3ex]
\ba{l}
q(A,z)\psi(u_{c-1}\otimes a)=u_c\cdot(\disc(x)a)\\[1ex]
q(\partial_A,z)\psi(u_{c-1}\otimes a)=u_c\cdot(\disc(y)a)
\ea
\quad\text{for $a\in e_{\det}H_ce$.}
\ea
\eneq
\item \label{rem3}
The diagrams \eqref{dia:Lc} and \eqref{dia:Lc-1} commute on $\X$.
\enum
\enremark
By Propositions \ref{prop:phiext}, \ref{prop:psiext}
and Remark \ref{rem:4.3} (iii),
we obtain the following proposition (see Proposition \ref{cor:A}).
\Prop\label{prop:fun}
Assume Condition \eqref{eq:HdiscH} holds.
Then we have isomorphisms 
of twisted $G$-equivariant $\W_\X$-modules with twist $c\tr:$
$$\vphi\cl\L_c\otimes_{eH_ce}eH_ce_{\det}
\iso\L_{c-1}\otimes \det(V)$$
and
$$\psi\cl\bl(\L_{c-1}\otimes \det(V)\br) \tensor_{e_{\det}H_ce_{\det}}
e_{\det}H_ce\iso\L_c.$$
\enprop

\subsubsection{}

Let us consider
$$\A_c=\bigl(p_*(\Endm_{\W_\X}(\L_c))^G\bigr)^\opp. $$
It is a \Ws-algebra on $\Hilb$ by Proposition \ref{prop:symred}.
Let $\A_c(0)$ be the subring of sections of order at most $0$.
For $m\in\Z$, $\L_{c+m}\otimes \det(V)^{\otimes-m}$ 
belongs to $\Mod_{c\tr}^G(\W_\X)$ (cf.\ \eqref{eq:lachi}).
Set
$$\A_{c,c+m}=(p_*\hom_{\W_\X}(\L_c,\L_{c+m}\otimes \det(V)^{\otimes-m}))^G.$$
Then $\A_{c,c+m}$ is an $(\A_c,\A_{c+m})$-bimodule.
Let
$\A_{c,c+m}(0)=(p_*\hom_{\W_\X(0)}(\L_c(0),
\L_{c+m}(0)\otimes \det(V)^{\otimes-m}))^G$.
Then $\A_{c,c+m}(0)$ is an $\A_c(0)$-lattice of
$\A_{c,c+m}$ and $\A_{c,c+m}(0)/\A_{c,c+m}(-1)\simeq L^{\otimes -m}$,
the associated line bundle on $\Hilb$
to $\CO_{\mu_\X^{-1}(0)}\otimes\det(V)^{\otimes -m}$
(cf.\ Proposition \ref{prop:symred} (iii)).

\subsection{Affinity of $\A_c$}
\label{secaff}
\subsubsection{}

As an application of Theorem \ref{gth:van},
we obtain the following vanishing theorem.
\Th\label{th:van}
Assume Condition \eqref{eq:HdiscH} holds
for $c+m$ \ro for all $m\in\Z_{>0}$\rf.
\bnum
\item
For any good $\A_c$-module $\M$,
$\prolim[K] H^i(K,\M)=0$ for $i>0$.
Here $K$ ranges over compact subsets of $\Hilb$.
\item
Any good $\A_c$-module $\M$
is generated by global sections on any compact subset of $\Hilb$.
\enum
\enth
\sketch
By Proposition \ref{prop:fun}, for any $m>0$,
$\L_{c+m}$ is a direct summand of 
a direct sum of copies of $\L_{c+m-1}\otimes \det(V)$
and
$\L_{c+m-1}\otimes\det(V)$ is a direct summand of a direct sum of 
copies of $\L_{c+m}$
in the category $\Mod_{(c+m)\tr}^{G}(\W_\X)$.
Hence
$\L_{c+m}\otimes \det(V)^{\otimes-m}$ is a direct summand of 
a direct sum of copies of $\L_{c}$
and
$\L_{c}$ is a direct summand of a direct sum of 
copies of $\L_{c+m}\otimes \det(V)^{\otimes-m}$
in the category $\Mod_{c\tr}^{G}(\W_\X)$ for any $m>0$.
It follows that $\A_{c,c+m}$ is a direct summand of
a direct sum of copies of $\A_{c}$
and
$\A_{c}$ is a direct summand of a direct sum of
copies of $\A_{c,c+m}$ for any $m>0$.
Moreover $\A_{c,c+m}$ is a good $\A_c$-module whose symbol is
$L^{\otimes-m}$.

Theorem \ref{gth:van} now gives the conclusion.
\QED

\subsubsection{}
Let us give an F-action on $\W_\X$ by
$\Fr_t(A_{ij})=tA_{ij}$,
$\Fr_t(\partial_{A_{ij}})=t^{-1}\partial_{A_{ij}}$,
$\Fr_t(z_{i})=tz_{i}$,
$\Fr_t(\partial_{z_{i}})=t^{-1}\partial_{z_{i}}$
and $\Fr_t(\h)=t^2\h$
for $t\in \Gm=\C^\times$.
Since $B_{ij}=\sigma_0(\h\partial_{A_{ji}})$ and
$\zeta_i=\sigma_0(\h\partial_{z_i})$,
the corresponding action of $\Gm$ on $\X$ is
$T_t((A,B,z,\zeta))=(tA,tB,tz,t\zeta)$.
Its induced $\Gm$-action on $\Hilb$
coincides with the action
induced by the scalar $\Gm$-action on $\C^2$.
We define the F-action on $\L_c$ by $\Fr_t(u_c)=u_c$.

Note that
\eqn\End_{\Mod_F(\W_\X[\h^{1/2}])}(\W_\X[\h^{1/2}])&\simeq&
\End_{\Mod_F(\W_{T^*(\g\times V)}[\h^{1/2}])}(\W_{T^*(\g\times V)}[\h^{1/2}])
\\&\simeq&
\C[\h^{-1/2}A_{ij},\h^{1/2}\partial_{A_{ij}},
\h^{-1/2}z_{i}, \h^{1/2}\partial_{z_{i}}]
\simeq\D(\g\times V).
\eneqn

The F-action on $\W_\X$ is compatible with the $G$-action on $\W$, and
hence $\A_c$ is also a \Ws-algebra on $\Hilb$ with F-action
(cf.\ Proposition~\ref{prop:symred} (iv)).
We define the F-action 
on $\L_{c-1}\otimes\det(V)$ by
$\Fr_t(u_{c-1}\otimes l)=t^{-n}u_{c-1}\otimes l$.
Hence $\A_{c,c-1}$ has a structure of $\A_c$-module with F-action.

\subsubsection{}
The $\bl((e+e_{\det})H_c(e+e_{\det})\br)^\opp$-module structure on 
$\L_c\oplus\bl(\L_{c-1}\otimes\det(V)\br)$ gives
a ring homomorphism 
$$(e+e_{\det})H_c(e+e_{\det})\To[\alpha] 
\End_{\A_c}(\A_c\oplus\A_{c,c-1})^\opp.$$
Since it is not compatible with the F-action,
we shall modify $\alpha$.

Set $$\tA_c=\A_c[\h^{1/2}]\quad\text{and}\quad
\tA_{c,c-1}=\A_{c,c-1}[\h^{1/2}].$$
Let $H_c\To[\beta]\corps[\h^{1/2}]\otimes_\C H_c$
be the ring homomorphism given by
$x_i\mapsto \h^{-1/2}\otimes x_i$,
$y_i\mapsto \h^{1/2}\otimes y_i$,
$w\mapsto 1\otimes w$ ($w\in W$).
\Lemma\label{lem:modact}
The composition
$$\Phi\cl (e+e_{\det})H_c(e+e_{\det})\To[\beta] 
\corps[\h^{1/2}]\otimes_\C (e+e_{\det})H_c(e+e_{\det})\To[\alpha]
\End_{\tA_c}(\tA_c\oplus\tA_{c,c-1})^\opp$$
sends $(e+e_{\det})H_c(e+e_{\det})$ to 
$\End_{\Mod_F(\tA_c)}(\tA_c\oplus\tA_{c,c-1})^\opp$.
\enlemma
\sketch
First let us show that $\Phi$ sends $eH_ce$
to $\End_{\Mod_F(\tA_c)}(\tA_c)^\opp$.
For a homogeneous element $a\in\C[\ts]^W$ of degree $k$, 
$\Phi(ae) (u_c)=\h^{-k/2}\tilde{a}(A)u_c$, where
$\tilde{a}(A)$ is the element of $\C[\g]^G$ such that
$\tilde{a}\vert_{\ts}=a$.
Since $\tilde{a}(A)$ is also homogeneous of degree $k$,
$\h^{-k/2}\tilde{a}(A)$ is $\Fr$-invariant,
and $\Phi(ae)$ belongs to $\Mod_F(\tA_c)$.
On the other hand,
we have $\Phi(\ysq e)(u_c)=\h\Delta_\g u_c$ and
$\h\Delta_\g$ is $\Fr$-invariant. Hence $\Phi(\ysq e)$ 
belongs to $\Mod_F(\tA_c)$.
Since $eH_ce$ is generated by $\C[\ts]^We$ and $\ysq e$, we have
$\Phi(eH_ce)\subset \End_{\Mod_F(\tA_c)}(\tA_c)$.

Similarly, we have 
$\Phi(e_{\det}H_ce_{\det})\subset \End_{\Mod_F(\tA_c)}(\tA_{c,c-1})$.

Let us show that
$\Phi(e\disc(x))\in \Hom_{\Mod_F(\tA_c)}(\tA_c,\tA_{c,c-1})$.
This follows from
$\Phi(e\disc(x))(u_c)=\h^{-n(n-1)/4}q(A,z)u_{c-1}\otimes l$,
$\Fr_t(q(A,z))=t^{n+n(n-1)/2}q(A,z)$ and
$\Fr_t(u_{c-1}\otimes l)=t^{-n}u_{c-1}\otimes l$.

For $a\in e_{\det}H_ce$, let us show that
$\Phi(a)\cl \tA_{c,c-1}\to \tA_c$ belongs to  $\Mod_F(\tA_c)$.
Since $\Phi(ae\disc(x))$ belongs to $\Mod_F(\tA_c)$,
and since $\Phi(e\disc(x))\vert_{\{q(A,z)\not=0\}}$ is an isomorphism in 
the category $\Mod_F(\tA_c\vert_{\{q(A,z)\not=0\}})$, it follows that
$\Phi(a)\vert_{\{q(A,z)\not=0\}}$ is in 
$\Mod_F(\tA_c\vert_{\{q(A,z)\not=0\}})$.
Hence we conclude that $\Phi(a)$ is in $\Mod_F(\tA_c)$.
Similarly, one shows that $\Phi(eH_ce_{\det})$
is contained in $\Hom_{\Mod_F(\tA_c)}(\tA_c,\tA_{c,c-1})$.
\QED
In particular we obtain a morphism of algebras
$$eH_ce\to \End_{\Mod_F(\tA_c)}(\tA_c)^\opp.$$
We denote by $\tilde{\vphi}$ and $\tilde{\psi}$ the modified morphisms
in $\Mod_F(\tA_c)$ given in Lemma~\ref{lem:modact}:
\eqn
\tilde{\vphi}&\cl&\tL_c\otimes_{eH_ce}eH_ce_{\det}
\To \tL_{c-1}\otimes\det(V),\\
\tilde{\psi}&\cl&\bl(\tL_{c-1}\otimes\det(V)\br)
\otimes_{e_{\det}H_ce_{\det}}e_{\det}H_ce\To \tL_c.
\eneqn

We define 
the order filtration $F(eH_ce)$ on $eH_ce$
by assigning order $1/2$ to $x_i$ and $y_i$.
Then the morphism
$eH_ce\to\End_{\Mod_F(\tA_c)}(\tA_c)^\opp$ 
is compatible with the order filtrations,
and the symbol map 
$\C[\ts\times\ts^*]^{W}\simeq\Gr^F(eH_ce)
\to \Gr^F\End_{\Mod_F(\tA_c)}(\tA_c)
\subset\Gamma(\Hilb,\CO_\Hilb)[\h^{\pm1/2}]$
coincides
with $\C[\ts\times\ts^*]^{W}_k
\To[\h^{-k}i_s]\h^{-k}\Gamma(\Hilb,\CO_\Hilb)$
by \eqref{eq:symb}.
Here $k\in\Z/2$.

\Lemma
The morphism
$eH_ce\to\End_{\Mod_F(\tA_c)}(\tA_c)^\opp$ is an isomorphism.
\enlemma
\sketch
Note that the subspace
$\Gr^F\End_{\Mod_F(\tA_c)}(\tA_c)\subset\Gamma(\Hilb,\CO_\Hilb)[\h^{\pm1/2}]$
is contained 
in $\oplus_{k\in\Z/2}\Gamma(\Hilb,\CO_\Hilb)_k\h^{-k}$
where $\Gamma(\Hilb,\CO_\Hilb)_k$ is the homogeneous part of weight $2k$
with respect to the $\Gm$-action.
Hence we have a chain of morphisms
\eqn
&&\C[\ts\times\Gt^*]^W\iso\Gr^F(eH_ce)\\
&&\hs{10ex}\to
\Gr^F(\End_{\Mod_F(\tA_c)}(\tA_c)^\opp)\hookrightarrow
\oplus_{k\in\Z/2}\Gamma(\Hilb,\CO_\Hilb)_k\h^{-k}
\iso\C[\ts\times\Gt^*]^W.\eneqn
Since the composition is the identity, the map
$\Gr^F(eH_ce)\to
\Gr^F(\End_{\Mod_F(\tA_c)}(\tA_c)^\opp)$
is bijective.
Hence the morphism
$eH_ce\to\End_{\Mod_F(\tA_c)}(\tA_c)^\opp$ is an isomorphism.
Note that $\bigcap_{k}F_k\bl(\End_{\Mod_F(\tA_c)}(\tA_c)\br)=0$.
\QED

\Remark
A similar argument shows that there is an isomorphism
$$eH_ce_{\det}\iso\Hom_{\Mod_F(\tA_c)}(\tA_c,\tA_{c,c-1}).$$
(See \S\;\ref{sec:quot}.)
\enremark

Let $\romano\in(\ts\times\ts^*)/W$ be the image of the origin of
$\ts\times\ts^*$. Then the Hilbert-Chow morphism
$\pi\cl \Hilb\to (\ts\times\ts^*)/W$ is $\C^\times$-equivariant, and
every point of $(\ts\times\ts^*)/W$ shrinks to $\romano$.

Now the following theorem is a consequence of Theorem \ref{th:Feq}.

\Th
Assume Condition \eqref{eq:HdiscH} holds
for $c+m$,  for all $m\in\Z_{>0}$ \ro this will be the case if
$c\not\in \frac{1}{n!}\BZ_{<0})$.
We have quasi-inverse equivalences
of categories between
$\Mod_F^\good(\tA_c)$ and $\Mod_{\coh}(eH_ce)$
\eqn
\Mod_F^\good(\tA_c)& \isotf & \Mod_{\coh}(eH_ce) \\
\M&\mapsto& \Hom_{\Mod_F^\good(\tA_c)}(\tA_c,\M)\\
\tA_c\otimes_{eH_ce} M & \mapsfrom & M.
\eneqn
\enth

Under this equivalence, $\tA_c$ and $\tA_{c,c-1}$
correspond to $eH_ce$ and $eH_ce_{\det}$, respectively.

\Th
Assume Condition \eqref{eq:HdiscH} holds
for $c+m$ \ro for all $m\in\Z_{>0}$\rf.
Assume also Condition \eqref{eq:HeH} holds
\ro these assumptions will be satisfied if
$c\not\in \frac{1}{n!}\BZ_{<0})$.
Let $\B_c=\Endomo_{\tA_c}(\tA_c\otimes_{eH_ce}eH_c)^\opp$.
We have quasi-inverse equivalences
of categories between
$\Mod_F^\good(\B_c)$ and $\Mod_{\coh}(H_c)$
\eqn
\Mod_F^\good(\B_c)& \isotf & \Mod_{\coh}(H_c) \\
\M&\mapsto& \Hom_{\Mod_F^\good(\B_c)}(\B_c,\M)\\
\B_c\otimes_{H_c} M & \mapsfrom & M.
\eneqn
\enth

\Remark
It would be very interesting to have a more direct construction of
$\tA_c\otimes_{eH_ce}eH_c$.
\enremark

\subsection{\Ws-algebras as fractions of $eH_ce$}
\label{sec:quot}
We explain how sections of $\tA_c$ over open subsets of $\Hilb$
can be obtained by inverting elements in the Cherednik algebra.

\smallskip
Let $\{F_j(H_c)\}_{j\in\Z/2}$ be the filtration of
$H_c$ consisting of elements of order $\le j$,
where we give order $1/2$ to $x_i$, $y_i$ and 
order $0$ to $w\in W$.
Then we have a canonical isomorphism
$\sigma\cl\Gr^F(H_c)\iso\C[\ts\times \ts^*]\rtimes W$.
We have induced filtrations on $eH_ce$ and $eH_ce_{\det}$,
and $\sigma$ induces isomorphisms
\eqn
\Gr^F(eH_ce)&\iso&\C[\ts\times \ts^*]^W,\\
\Gr^F(eH_ce_{\det})&\iso&\C[\ts\times \ts^*]^{W,\,\det}.
\eneqn
Composing with the morphism
$\C[\ts\times \ts^*]\to \C[\ts\times \ts^*][\h^{-1/2}]$
given by $a(x,y)\mapsto a(\h^{-1/2}x,\h^{-1/2}y)$,
we obtain homomorphisms
\eqn
\Gr^F(eH_ce)&\To&\C[\ts\times \ts^*]^W[\h^{-1/2}],\\
\Gr^F(eH_ce_{\det})&\To&\C[\ts\times \ts^*]^{W,\,\det}[\h^{-1/2}].
\eneqn

We shall set $\tW_\X=\W_\X[\h^{1/2}]$
and $\tW_\X(0)=\W_\X(0)+\h^{1/2}\W_\X(0)$.
We set $\tL_c=\tW_\X\otimes_{\W_\X}\L_c$.
Then $\tL_c\oplus\tL_{c-1}\otimes\det(V)$ has a structure
of $(\tW_\X,(e+e_{\det})H_c(e+e_{\det}))$-bimodule.
The action of $eH_{c}e_{\det}$ is given by
$\tilde{\vphi}\cl 
\tL_{c}\otimes_{eH_{c}e}eH_{c}e_{\det}\to\tL_{c-1}\otimes \det(V)$.
On the other hand, we have canonical isomorphisms
$\Gr^F(\tL_{c})\simeq\Gr^F(\tL_{c-1})\iso\CO_{\mu_\X^{-1}(0)}[\h^{\pm1/2}]$.
Here $F(\tL_c)$ (resp.\ $F(\tL_{c-1})$)
is the order filtration given by
$F_k(\tL_c)=\h^{-k}\tW_\X(0)u_c$ 
(resp.\ $F_k(\tL_{c-1})=\h^{-k}\tW_\X(0)u_{c-1}$) for $k\in\Z/2$.

We have a commutative diagram:
\eq&&\ba{c}
\xymatrix@R=3ex@C=7ex{
{\Gr^F(\tL_{c})\otimes\Gr^F(eH_{c}e)}\ar[r]\ar[d]
&{\CO_{\mu_\X^{-1}(0)}[\h^{\pm1/2}]\otimes \C[\ts\times\ts^*]^{W}[\h^{-1/2}]}
\ar[dd]^{i_s}\\
{\Gr^F(\tL_{c}\otimes eH_{c}e)}\ar[d]
\\
{\Gr^F(\tL_{c})}\ar[r]^\sim&{\CO_{\mu_\X^{-1}(0)}[\h^{\pm1/2}].}
}\ea
\eneq

The morphism $\tilde\vphi$ is order-preserving 
and we obtain a commutative diagram
\eq&&\ba{c}
\xymatrix@R=3ex@C=7ex{
{\Gr^F(\tL_{c})\otimes\Gr^F(eH_{c}e_{\det})}\ar[r]\ar[d]&
{\CO_{\mu_\X^{-1}(0)}[\h^{\pm1/2}]
\otimes\C[\ts\times\ts^*]^{W,\,\det}[\h^{-1/2}]}
\ar[dd]^{i_d}\\
{\Gr^F(\tL_{c}\otimes eH_{c}e_{\det})}\ar[d]^{\tilde{\vphi}}\\
{\Gr^F(\tL_{c-1}\otimes\det(V))}\ar[r]^\sim
&{\CO_{\mu_\X^{-1}(0)}[\h^{\pm1/2}]\otimes\det(V).}
}\ea
\eneq
Hence, for any $a\in eH_{c}e_{\det}$, the morphism
$a\cl \tL_{c}\to\tL_{c-1}\otimes\det(V)$
is an isomorphism on $\{i_d(\sigma(a))\not=0\}$.
Then, for $b\in  eH_{c}e_{\det}$, we can define
$$ba^{-1}\in\End_{\Mod^G_{F,\,c\tr}(\tW_X\mid_{\{i_d(\sigma(a))\not=0\}})}
(\tL_{c}\mid_{\{i_d(\sigma(a))\not=0\}})^\opp$$
as the composition
$$\xymatrix@R=.5ex@C=9ex{
{\tL_c}\ar[rd]^(.3)b\ar@{.>}[dd]_{ba^{-1}}\\
&{\tL_{c-1}\otimes\det(V).}\\
{\tL_c}\ar[ru]_(.3)a
}
$$
Thus we obtain
$ba^{-1}$ as an F-invariant section of
$\tA_c$ defined on $\{i_d(\sigma(a))\not=0\}$.
Note that $ba^{-1}=bc(ac)^{-1}$
for a non-zero element $c\in e_{\det}H_ce$.
Note also that the image of $ac\in eH_ce$ in $\Gamma(\Hilb;\tA_c)$
is invertible only on $\{i_d(\sigma(a))\not=0\}\cap\{i_d(\sigma(c))\not=0\}
\cap(\Hilb\setminus E)$.

\Remark
The morphism $\tilde{\psi}\cl\bl(\tL_{c-1}\otimes\det(V)\br)
\otimes_{e_{\det}H_ce_{\det}}e_{\det}H_ce\to\tL_c$ is also order-preserving,
and it induces a commutative diagram
\eqn
\xymatrix@R=3ex@C=3ex{
{\Gr^F(\tL_{c-1}\otimes\det(V))\otimes \Gr^F(e_{\det}H_{c}e)}\ar[r]\ar[d]
&{\CO_{\mu_\X^{-1}(0)}[\h^{\pm1/2}]\otimes\det(V)
\otimes\C[\ts\times\ts^*]^{W,\,\det}[\h^{-1/2}]}\ar[dd]^-{\tau\cdot i_d}\\
{\Gr^F(\tL_{c-1}\otimes\det(V)\otimes e_{\det}H_{c}e)}\ar[d]^-{\tilde{\psi}}
\\
{\Gr^F(\tL_{c})}\ar[r]^\sim&{\CO_{\mu_\X^{-1}(0)}[\h^{\pm1/2}].}
}
\eneqn
Hence, for any $b\in  e_{\det}H_{c}e$,
the morphism $b\cl \tL_{c-1}\otimes\det(V)\to\tL_c$
is never an isomorphism on the exceptional divisor $E$.
\enremark

\subsection{Rank $2$ case}
\label{sec:rank2}
Let us consider the case $n=2$. 
Let $x_0=x_1+x_2$, $x=x_1-x_2$,
$y_0=(y_1+y_2)/2$ and $y=(y_1-y_2)/2\in H_c$.
Then $[y_0,x_0]=1$, $[y,x]=1-2cs$ where
$s=s_{12}$.
Since $y$, $x$ and $s$ commute with $\C[x_0,y_0]$,
we have an isomorphism of algebras
$\C[x_0,y_0]\otimes H'_c\iso H_c$, where
$H'_c$ is the subalgebra of $H_c$ generated by $x$, $y$ and $s$.

We have
\eqn
&&eH_ce_{\det}H_ce=eH_ce\Longleftrightarrow
H_ce_{\det}H_c=H_c\Longleftrightarrow c\not=1/2,\\
&&e_{\det}H_ceH_ce_{\det}=e_{\det}H_ce_{\det}\Longleftrightarrow
H_ceH_c=H_c\Longleftrightarrow c\not=-1/2.
\eneqn
Indeed, the first equivalences follow from
the fact that $ye_{\det}x-xe_{\det}y=e[y,x]=(1-2c)e$ and when $c=1/2$,
there is a one-dimensional representation with
$x,y\mapsto0$, $s\mapsto1$.
The second follows from the first by the isomorphism
$H_c\simeq H_{-c}$ given by $s\mapsto -s$.
It follows that Condition \eqref{eq:HdiscH}
is satisfied for all $c+n$ ($n\in\Z_{>0}$) 
if and only if $c\not=-1/2,-3/2,\ldots\,$.

\medskip
Note that $x$, $y\in \C[\ts\times\ts^*]^{W,\det}$ and
$\Hilb=\{i_d(x)\not=0\}\cup\{i_d(y)\not=0\}$,
because $\mu_\X^{-1}(0)\cap\{q(A,z)=q(B,z)=0\}\subset
\set{(A,B,z,0)\in\X}{Az,Bz\in\C z}=\emptyset$.
Quantized symplectic coordinates of $\tA_c$
are given by
$$ ((ey)(ex)^{-1},\h^{1/2}ex_0;-\hbar ex^2/2,\h^{1/2} ey_0)\quad\text{on }
\{i_d(x)\not=0\}$$
 and
$$((ex)(ey)^{-1},\h^{1/2}ex_0;\hbar ey^2/2,\h^{1/2} ey_0)
\quad\text{on } \{i_d(y)\not=0\}.$$
Indeed, we have $[-ex^2/2,(ey)(ex)^{-1}]=e$,
because
$$(ey)(ex)^{-1}(ex^2)=(ey)(ex)^{-1}(ex)(e_{\det}x)=eyx\quad \text{and}$$
$$(ex^2)(ey)(ex)^{-1}=(ex^2y)(ex)^{-1}=(eyx^2-2ex)(ex)^{-1}=
(eyx)(ex)(ex)^{-1}-2e=eyx-2e.$$

Note that this provides an isomorphism
$\Hilb\iso T^*(\Pbb^1\times \C)$.
The projection
$\Hilb\to\Pbb^1$ is given by $[i_d(x):i_d(y)]$
with the notation of homogeneous coordinates.
By the isomorphism above, we have $E\simeq T^*_{\Pbb^1}\Pbb^1\times T^*\C$.

Note that $(xe)^{-1}(ye)$ is invertible only on
$\{i_s(x^2)\not=0\}=\{i_d(x)\not=0\}\setminus E$ for $c\not=-1/2$, because
$exyx=ex(xy+1-2cs)=ex^2y+(1+2c)ex$ and
$(xe)^{-1}(ye)=(x^2e)^{-1}(xye)=(ex^2)^{-1}(exyx)(ex)^{-1}
=(ey)(ex)^{-1}+(1+2c)(ex^2)^{-1}$.

\medskip
Set $(a,\partial_a)=((ey)(ex)^{-1},-ex^2/2)$ and
$(b,\partial_b)=(((ex)(ey)^{-1},ey^2/2)$ and $\la=c-1/2$.
Then we have
\eq
b=a^{-1}\quad\text{and}\quad\partial_b=-a(a\partial_a-\la).
\label{eq:patchA}
\eneq
Indeed, we have
\eqn
-a(a\partial_a-\la)&=&(ey)(ex)^{-1}\bl((ey)(ex)^{-1}(ex^2)/2+c-1/2\br)\\
&=&(1/2)(ey)(ex)^{-1}(eyx+2c-1)
=(1/2)(ey)(ex)^{-1}(exy)=ey^2/2.
\eneqn

Recall that $\romano\in(\ts\times\Gt^*)/W$
is the image of the origin of $\ts\times\Gt^*$.
The inverse image $\pi^{-1}(\romano)$
by the Hilbert-Chow morphism $\pi$
is $T^*_{\Pbb^1}\Pbb^1\times\{0\}\subset T^*\Pbb^1\times T^*\C$.
We identify it with $\Pbb^1$.
Then, \eqref{eq:patchA} gives an isomorphism
$$\Endomo_F(\tA_c)\mid_{\pi^{-1}(\romano)}\iso\D_{\Pbb^1,\,\lambda}
\otimes \C[x_0,y_0]$$
with $\la=c-1/2$.
Here, $\D_{\Pbb^1,\,\lambda}$ is the twisted ring of differential operators
(e.g.\ see \cite[\S\,2]{K}). If $\lambda$ is an integer, then
$\D_{\Pbb^1,\,\lambda}\simeq\CO_{\Pbb^1}(\la)\otimes \D_{\Pbb^1}\otimes
\CO_{\Pbb^1}(-\la)$.
Hence, we have 
a ring isomorphism $eH'_ce\simeq\Gamma(\Pbb^1;\D_{\Pbb^1,\,\lambda})$
and an equivalence
$\Mod^\good_F(\tA_c)\simeq 
\Mod_\good(\D_{\Pbb^1,\,\lambda}\otimes \C[x_0,y_0])$.
It is well-known (cf.\ e.g. \cite[\S\,7]{K}) that
$\Mod_\good(\D_{\Pbb^1,\,\lambda})$ is equivalent to
$\Mod_\coh(\Gamma(\Pbb^1;\D_{\Pbb^1,\,\lambda}))$
if and only if $\la\not=-1,-2,\ldots$ (i.e.\ $c\not=-1/2,-3/2,\ldots$).

\end{document}